\documentclass[12pt,reqno]{article}
\usepackage{amsfonts,mathrsfs}
\usepackage{fancyhdr,amscd}
\usepackage{anysize}
\usepackage{setspace}
\usepackage{graphicx,epsfig}
\usepackage{amsthm,latexsym,amsmath,amssymb}
\usepackage[perpage,symbol*]{footmisc}
\usepackage[colorlinks=true]{hyperref}
\usepackage{xcolor}
\usepackage{amsfonts,amsmath,amssymb,latexsym}
\usepackage{mathrsfs}

\usepackage{epsfig,CJK,fancyhdr}
\usepackage{CJK,fancyhdr,amscd}

\newtheorem{thm}{Theorem}[section]
\newtheorem{defi}[thm]{Definition}
\newtheorem{prop}[thm]{Proposition}
\newtheorem{lem}[thm]{Lemma}

\newtheorem{rem}[thm]{Remark}

\makeatletter \@addtoreset{equation}{section} \makeatother
\def\beq{\begin{equation}}
\def\eeq{\end{equation}}
\def\endproof{$\hfill\Box$}

\begin{document}
\baselineskip=20pt  \hoffset=-3cm \voffset=0cm \oddsidemargin=3.2cm
\evensidemargin=3.2cm \thispagestyle{empty}\vspace{10cm}

\hbadness=10000
\tolerance=10000
\hfuzz=150pt
\title{\textbf{An index theory with applications to homoclinic  Orbits of Hamiltonian systems  and Dirac equations}}
\author{\Large Qi Wang $^{{\rm a}}$,$\quad$  Chungen Liu $^{{\rm b}}$}
\date{} \maketitle
\begin{center}
\it\scriptsize ${}^{\rm a}$Institute of Contemporary Mathematics, School of Mathematics and Information Science,\\ Henan
University, Kaifeng 475000, PR China\\
${}^{\rm b}$School of Mathematics and Information Science, Guangzhou University, Guangzhou 510006, PR China
\end{center}

\footnotetext[0]{$^a$Partially supported by NNSF of China(11301148).}
\footnotetext[0]{$^{\rm b}${\bf Corresponding author.} Partially  supported by NNSF of China(11471170).}
\footnotetext[0]{\footnotesize\;{\it E-mail address}: mathwq@henu.edu.cn. (Qi Wang), liucg@nankai.edu.cn (Chungen Liu).}

\date{} \maketitle

\noindent
{\bf Abstract:} {\small In this paper, we will define the index pair $(i_A(B),\nu_A(B))$ by the  dual variational method, and show  the relationship   between the  indices defined by different methods. As applications, we  apply the index $(i_A(B),\nu_A(B))$ to study the existence and multiplicity of homoclinic orbits of nonlinear Hamiltonian systems and solutions of nonlinear Dirac equations.}


\section{Introduction}\label{section-introduction}

Many problems can be displayed as a self-adjoint operator equation
\[
 Au=F'(u),\;u\in D(A)\subset \mathbf H,\eqno{(OE)}
\]
where $\mathbf H$ is an infinite-dimensional separable Hilbert space, $A$ is a  self-adjoint operator on $\mathbf H$ with its domain $D(A)$, $F$ is a nonlinear functional on $\mathbf H$, such as Dirichlet problem
for Laplace's equation on bounded domain, periodic solutions of Hamiltonian systems, nonlinear Dirac equations, system of diffusion equations, Schr\"{o}dinger equation, periodic solutions of wave equation and so on. By variational
method, we know that the solutions of (OE) correspond to the critical points of a functional on a Hilbert space. So we can transform the problem of finding the solutions of
(OE) into the problem of finding the critical points of a functional.  Many theories have been developed to do so. Among these theories, Morse theory is one of the remarkable theories, and it has a great advantage in displaying the relationship between the global and local behavior of the functional.

Morse theory can be used directly in Dirichlet problem for Laplace's equation on bounded domain and periodic solutions of second order Hamiltonian systems,
since the Morse indices of the critical points are finite. But for the problems of periodic solutions of first order Hamiltonian systems,  Schr\"{o}dinger equations, wave equations,
Morse theory cannot be used directly because  in these situations the functionals are strongly indefinite in the sense that they are  unbounded from above and below and the Morse indices at the critical points of these functionals are infinite.
 Fortunately, some methods have been developed to deal with these situations, such as Galerkin approximation methods,
saddle point reduction(a kind of Lyapunov-Schmidt procedure, see e.g Amann\cite{Amann-1976}, Amann and Zehnder\cite{Amann-Zehnder-1980} and Chang\cite{Chang-1993}), dual variational methods and convex analysis theory (see e.g Aubin and Ekeland \cite{Aubin-Ekeland-1984}, Ekeland\cite{Ekeland-1990},
Ekeland and Temam \cite{Ekeland-Temam-1976}). By these methods, the solutions of (OE) correspond to the critical points of functionals with finite relative Morse indices, then one can use Morse theory to find the solutions of (OE).

 Related to Morse theory, the relative index theory is worth to pay close attention.  By the work  \cite{Ekeland-1984} of  I. Ekeland, an index theory for convex linear Hamiltonian systems was established.
By the works  \cite{Conley-Zehnder-1984,Long-1990,Long-1997,Long-Zehnder-1990} of Conley, Zehnder and Long, an index theory for symplectic paths was introduced.
These index theories have important and extensive applications,
e.g \cite{Dong-Long-1997,Ekeland-Hofer-1985,Ekeland-Hofer-1987,Liu-Long-Zhu-2002,Long-Zhu-2000}.
In \cite{Long-Zhu-1999, Long-Zhu-2000-2} Long and Zhu defined spectral flows for paths of linear operators and redefined Maslov index for symplectic paths.
Additionally,  Abbondandolo defined a relative Morse index theory for Fredholm operator with compact perturbation (see\cite{Abb-2001} and the references therein). Chen and Hu defined the Maslov index for homoclinic orbits of Hamiltonian systems in \cite{Chen-Hu-2007}.
In the study of the $L$-solutions (the solutions starting and ending at the same Lagrangian subspace $L$) of Hamiltonian systems, the second author
of this paper introduced in  \cite{Liu-2007} an index theory for symplectic paths using the algebraic methods and gave some applications in
 \cite{ Liu-2007, Liu-2007-2}. Then this index had been generalized by the authors of this paper and Lin in  \cite{Liu-Wang-Lin-2011}.

In addition to the above index theories defined for specific forms, Dong in \cite{Dong-2010} developed an index theory for abstract operator equations (OE).
As an essential condition, he assumed that the embedding $D(A)\hookrightarrow \mathbf H$ was compact.
As applications, he considered the second order Hamiltonian systems, elliptic partial differential equations and first order Hamiltonian systems.
Recently, the authors of this paper in \cite{Wang-Liu-2014, Wang-Liu-2015} defined their index theory for abstract operator equations (OE)
by relative Fredholm index and spectral flow. We also needed the condition of  compact embedding $D(A)\hookrightarrow \mathbf H$.
As applications, we considered delay differential system and a kind of infinite dimensional Hamiltonian systems.
But for the cases of wave equations, beam equations and so on, the above two index theories will not work,
since the corresponding operators have essential spectrum and the condition of compact embedding will not be satisfied.

In order to overcome this difficulty, if the operator $A$ has no compact resolvent, in \cite{Wang-Liu-2016} we considered three cases in this situation, and defined the index pairs $(i^+_A(B),\nu^+_A(B))$, $(i^-_A(B),\nu^-_A(B))$ and $(i^0_A(B),\nu^0_A(B))$. Roughly speaking, we defined the index pairs $(i^+_A(B),\nu^+_A(B))$ and  $(i^-_A(B),\nu^-_A(B))$ by the method of dual variational and defined the index pair $(i^0_A(B),\nu^0_A(B))$ by the method of saddle point reduction. We  give the applications of $(i^+_A(B),\nu^+_A(B))$ and  $(i^-_A(B),\nu^-_A(B))$ for wave equation but have no applications of $(i^0_A(B),\nu^0_A(B))$.

 We now return  to the third case considered in \cite{Wang-Liu-2016} where the index pair $(i^0_A(B),\nu^0_A(B))$ was defined via saddle point reduction method. In this paper, we will define an index pair $(i_A(B),\nu_A(B))$ by using the dual variational method just for the operator pair $(A,B)$  (see Section \ref{section-the definition of index pair}).
Briefly speaking, let $\mathbf{H}$  be an infinite dimensional separable Hilbert space, $A$ is an unbounded self-adjoint operator on $\mathbf H$ with its essential spectrum $\sigma_e(A)\cap(a,b)=\emptyset$,
 for any bounded self-adjoint operator $B$ on $\mathbf H$ with its spectrum satisfying $\sigma(B)\in (a,b)$, we will define the index pair $(i_A(B),\nu_A(B))$.
 Of course, we will show  the relationship  between the indices defined by different methods and give  the relation of the  index theory  with the spectral flow of the related operator. Finally, we  apply the index $(i_A(B),\nu_A(B))$ to study the existence and multiplicity of homoclinic orbits of nonlinear Hamiltonian systems and solutions of nonlinear Dirac equations in Section \ref{section-Proof of main results}. \\
  \noindent{\textit{ Application A:  existence and multiplicity of  homoclinic orbits of nonlinear Hamiltonian system.}}$\;$ Consider the first order Hamiltonian system
\[
\left\{\begin{array}{ll}
 \dot{z}(t)=\mathcal{J}\nabla_zH(t,z),\\
\displaystyle\lim_{t\to\infty}z(t)=0,
  \end{array}
\right.\eqno{(HS)}
\]
where $z\in\mathbb{R}^{2N}$, $\mathcal{J}=\left(\begin{matrix}
              0 &-I_N\\
             I_N &0
             \end{matrix}\right)$ with $I_N$ the identity map on $\mathbb{R}^N$ and $H\in C^1(\mathbb{R}\times\mathbb{R}^{2N},\mathbb{R})$.
The solutions of (HS) are called homoclinic orbits of nonlinear Hamiltonian system. As a special case of dynamical systems, Hamiltonian systems are very important in the study of gas dynamics, fluid mechanics, relativistic mechanics and nuclear physics. However it is well known that homoclinic solutions play an important role in analyzing the chaos of Hamiltonian systems. If a
system has the transversely intersected homoclinic solutions, then it must be chaotic. If it has  smooth  connected homoclinic solutions, then it cannot stand the perturbation, and its perturbed system probably produces chaotic phenomena. Therefore, it is of practical importance and mathematical significance to consider the existence of homoclinic solutions of Hamiltonian systems emanating from the origin. The existence and multiplicity of homoclinic orbits for the first order system were studied extensively by means of critical point theory, and many results were obtained
under the assumption that $H(t,z)$ depends periodically on $t$ (see \cite{Arioli-Szulkin-1999,Chen-Ma-2011,Ding-2006,Ding-Girardi-1999,Ding-Willem-1999,Hofer-Wysocki-1990,S-1992,S-1993,Sun-Chu-Feng-2013,Szulkin-Zou-2001,Tanaka-1991,Wang-Xu-Zhang-2010,Wang-Xu-Zhang-2012,Zelati-Ekeland-1990} and the references therein). As authors known, the periodicity is used to protect some kind of compactness such as the (PS) condition.  Without assumptions of periodicity the problem is quite different in nature.
 To the best of our knowledge, the authors in \cite{Ding-Li-1995} firstly obtained the existence of homoclinic orbits for a class of first order systems in the non-periodic case. They assume the Hamiltonian function $H(t,z)$ has the following form
 \begin{equation}\label{eq-the form of Hamiltonial function H}
 H(t,z)=\frac{1}{2}(L(t)z,z)+R(t,z).
 \end{equation} provided that $L$ has a special form and $R(t,z)$ satisfies some kind of superquadratic or subquadratic growth conditions at infinity with respect to $z$. Then in \cite{Ding-Jeanjean-2007} by assuming that $L$ satisfies a more general condition and  $R(t, z)$ is asymptotically quadratic at infinity with respect to $z$, they also obtained the existence and multiplicity of homoclinic orbits. Additionally, there are
also a few papers devoted to the non-periodic case(see \cite{Ding-Lee-2009,Sun-Chen-Nieto-2011, Zhang-Tang-Zhang-2013,Zhang-Tang-Zhang-2015}).

In \cite{Ding-Jeanjean-2007}, for any $2N\times 2N$ matrix $M$, they say $ M\geq 0$ if  and only if
\[
 \displaystyle\min_{\xi\in\mathbb{R}^{2N},|\xi|=1}M\xi\cdot \xi\geq 0,
\]
 and denote $M\ngeq 0$ if and only if $M\geq 0$ does ont hold. Assume \\
(L). There exists $b>0$ such that the set $\Lambda^b:=\{t\in\mathbb{R}:\mathcal{J}_0L(t)-b\ngeq 0\}$ is nonempty and has finite measure, where $\mathcal{J}_0=\left(\begin{matrix}
              0 &I_N\\
             I_N &0
             \end{matrix}\right)$.\\
Assume $H$ has the form \eqref{eq-the form of Hamiltonial function H}, let
\begin{equation}\label{eq-the definition of operator A in HS}
A:=-(\mathcal{J}\frac{d}{dt}+L(t)),
\end{equation}
 then $A$ is self-adjoint on $L^2(\mathbb{R},\mathbb{R}^{2N})$. We have the following lemma.
\begin{lem}\label{Proposition 2.1-Ding-Jeanjean-2007}\cite[Proposition 2.1]{Ding-Jeanjean-2007}
 Assume (L) is satisfied, then
\[
 \sigma_e(A)\subset \mathbb{R}\setminus (-b_{max},b_{max})
\]
 with
\[
 b_{max}:=\sup\{b:|\Lambda^b|<\infty\}.
\]
\end{lem}
 In this part let the Hilbert space $\mathbf H:=L^2(\mathbb{R},\mathbb{R}^{2N})$, $\mathcal{L}_s(\mathbf H)$ the bounded self-adjoint operators on $\mathbf H$ and
 \[
 \mathcal{L}_s(\mathbf H,-b_{max},b_{max}):=\{B\in \mathcal{L}_s(\mathbf H)|\sigma(B)\in (-b_{max},b_{max})\}.
 \]
  Lemma \ref{Proposition 2.1-Ding-Jeanjean-2007} motivates us that our index $(i_A(B),\nu_A(B))\;(B\in \mathcal{L}_s(\mathbf H,-b_{max},b_{max}))$ can be used here to study the existence and multiplicity of homoclinic orbits.
 Denote $\mathcal{L}_s(\mathbb{R}^{2N})$ the set of all $2N\times 2N$ symmetric matrices and   $\mathcal{B}\subset C(\mathbb{R},\mathcal{L}_s(\mathbb{R}^{2N}))$ the set of all bounded
 symmetric $2N\times 2N$ matrix functions.
 For any $B\in\mathcal{B}$, it is easy to see $B$ determines a bounded self-adjoint operator on $\mathbf H$, by
 \begin{equation}\label{eq-multiple operator g}
  z(t)\mapsto B(t)z(t), \;\forall z\in \mathbf H,
 \end{equation}
 without confusion, we still denote this operator by $B$, that is to say
 we have the continuous embedding $\mathcal{B}\hookrightarrow \mathcal{L}_s(\mathbf H)$.
Assume  $R\in C^2(\mathbb{R}\times\mathbb{R}^{2N},\mathbb{R})$ in (\ref{eq-the form of Hamiltonial function H}) satisfying the following conditions
\begin{list}{}
{\setlength{\topsep}{1ex} \setlength{\itemsep}{0ex}
\setlength{\leftmargin}{2.6em}
 }
\item[$(R_0)$] $\nabla_zR(t,0)\equiv 0$, and  $B_0{ :=}\nabla^2_zR(t,0)\in  \mathcal{B}\cap\mathcal{L}_s(\mathbf H,-b_{max},b_{max})$.

 \item[$(R_1)$] There exists a constant $ \delta>0$ such that
\[
 (-b_{max}+\delta)Id<\nabla^2_zR(t,z)<(b_{max}-\delta)Id,\;\forall  (t,z)\in \mathbb{R}\times\mathbb{R}^{2N}.
\]

\item[$(R^\pm_2)$] There exist $B_1,B_2\in \mathcal{B}\cap \mathcal{L}_s(\mathbf H,-b_{max},b_{max})$ with $\pm B_1< \pm B_2$,
$i_A(B_1)=i_A(B_2)$ and $\nu_A(B_1)=\nu_A(B_2)=0$,
such that
\[
\pm \nabla^2_zR(t,z)\geq  \pm B_1(t),\;\forall (t,z)\in \mathbb{R}\times\mathbb{R}^{2N}
\]
and
\[
\pm \nabla^2_zR(t,z)\leq \pm B_2(t),\; |(t,z)|>K
\]
 for some constant $K>0$.
\end{list}
Then we have the following result.
\begin{thm}\label{thm-asymptotically linear result 1}
 Assume (L) is satisfied,  $R\in C^2(\mathbb{R}\times\mathbb{R}^{2N},\mathbb{R})$ satisfies conditions $(R_0)$, $(R_1)$ and $(R^\pm_2)$, if
\[
 i_A(B_0)>i_A(B_2)(or\; i_A(B_0)+\nu(B_0)<i_A(B_1)),
\]
then (HS) has a nontrivial homoclinic orbit.  Further more, if $R$ satisfies more conditions, we will get more results.\\
\noindent  A. If $\nabla^2_zR$ is globally Lipschitz continuous on $z$, that is to say there exists a constant $L_R>0$,such that
\[
|\nabla^2_zR(t,z_1)-\nabla^2_zR(t,z_2)|\leq L_R|z_1-z_2|,\;\forall (t,z)\in\mathbb{R}\times\mathbb{R}^{2N},
\]
 and $\nu_A(B_0)=0$, then (HS) has another nontrivial homoclinic orbit different from the above one. \\
\noindent B. If $R$ is even in $z$, then (HS) has $i_A(B_0)-i_A(B_2)${\rm(}or $i_A(B_1)-i_A(B_0)-\nu(B_0)${\rm)} pairs of nontrivial homoclinic orbit.
\end{thm}

  Compared to the known results, we note that in   \cite{Ding-Jeanjean-2007,Sun-Chen-Nieto-2011}, where the condition $R\geq 0$ is required and the authors  used some spectral requirements to act as the twisting conditions at the origin and the infinity on $R$. Here we use the indices to state the  the twisting conditions. We note that in \cite{Chen-Hu-2007}, the authors developed an index theory for homoclinic solutions of  first order Hamiltonian systems but with no application to study the existence and multiplicity. Our this work is the first one to apply index theory to study the existence and multiplicity of homoclinic solutions for  first order Hamiltonian systems without compactness assumption.\\

\noindent {\textit{Application B:  existence and multiplicity of solutions of Nonlinear Dirac equations.}}$\;$
Nonlinear Dirac equations occur in the attempt to model extended relativistic particles with external fields, in a general form, such equations are given by
\[
 -i\hbar\partial_t\psi=ic\hbar\displaystyle\sum^3_{k=1}\alpha_k\partial_k\psi-mc^2\beta\psi-M(x)\psi+G_\psi(x,\psi),\eqno(D)
\]
where $x=(x_1,x_2,x_3)\in\mathbb{R}^3$, $\partial_k=\frac{\partial}{\partial_{x_k}}$, $c$ denotes the speed of light, $m>0$ is the mass of the electron, $\hbar$ denotes Planck's constant,
$M(x)$ is the matrix potential and in the nonlinearity term $G:\mathbb{R}^3\times \mathbb{C}^4\to\mathbb{R}$ represents a nonlinear self-coupling.
A solution $\psi:\mathbb{R}\times\mathbb{R}^3\to\mathbb{C}^4\ $ of (D) is a wave function which represents the state of a relativistic electron.
Furthermore, $\alpha_1,\alpha_2,\alpha_3$ and $\beta$ are $4\times 4$ complex matrices whose standard form (in $2\times 2$ blocks) is
\[
 \beta=\left(
\begin{matrix}
 I& 0\\
0&-I
\end{matrix}
\right),\;
\alpha_k=\left(
\begin{matrix}
 0& \sigma_k\\
\sigma_k&0
\end{matrix}
\right),\;k=1,2,3
\]
with
\[
 \sigma_1=\left(
\begin{matrix}
 0& 1\\
1&0
\end{matrix}
\right),\;
\sigma_2=\left(
\begin{matrix}
 0& -i\\
i&0
\end{matrix}
\right),\;
\sigma_3=\left(
\begin{matrix}
 1& 0\\
0&-1
\end{matrix}
\right).\;
\]
One verifies that $\beta=\beta^*$, $\alpha_k=\alpha^*_k$, $\alpha_k\alpha_l+\alpha_l\alpha_k=2\delta_{kl}$ and $\alpha_k\beta+\beta\alpha_k=0$,
due to these relations, the linear operator $\mathcal{H}_0:=-ic\hbar\displaystyle\sum^3_{k=1}\alpha_k\partial_k+mc^2\beta$
is a symmetric operator, such that
$
 \mathcal{H}_0^2=-c^2\hbar^2\Delta+m^2c^4.
$
The stationary solutions of equation ($D$) are found by the Ansatz
$
 \psi(t,x)=e^{\frac{i\theta t}\hbar}z(x).
$
Then $z:\mathbb{R}^3\to\mathbb{C}^4$ satisfies the equation
\[
 -ic\hbar\displaystyle\sum^3_{k=1}\alpha_k\partial_kz+mc^2\beta z+M(x)z=G_z(x,z)-\theta z.
\]
Now we re-written  the above equation by
\[
 -i\displaystyle\sum^3_{k=1}\alpha_k\partial_kz+V(x)\beta z=H_z(x,z).\eqno(DE)
\]
with $\alpha_*$($*=1,2,3$) and $\beta$ defined above, $V:\mathbb{R}^3\to\mathbb{R}$ and $H:\mathbb{R}^3\times\mathbb{C}^4\to \mathbb{R}$.

In \cite{Merle-1988}, Merle study the problem (DE) with a constant potential $V(x)=\omega$ and nonlinear term $F$ representing the so called Soler model. In \cite{Esteban-Sere-1995}, it seems that Esteban and S\'{e}r\'{e} were pioneers in using variational methods to study Soler model. But it's worth to note that the method used in \cite{Esteban-Sere-1995} doesn't work in the case of non-autonomous  systems which are important in quantum mechanics. Then Bartsch and Ding in \cite{Bartsch-Ding-2006} studied the existence and multiplicity of  the non-autonomous Dirac equations by their critical point theories for strongly indefinite problems. Compared with the periodic assumption of $V$ and $H$ in \cite{Bartsch-Ding-2006}, Ding and Ruf in \cite{Ding-Ruf-2008}  studied the existence and multiplicity of solutions of (DE) with non-periodic assumption.  After that, there are many works dedicated to study the Dirac equation with $V$ and  $H$ satisfying several different hypotheses (see \cite{Ding-Liu-2007,Ding-Xu-2015,Figueiredo-Pimenta-2017} and the references therein).

Motivated by \cite{Bartsch-Ding-2006} and \cite{Ding-Ruf-2008}, we study the existence and multiplicity of solutions of (DE) with non-periodic assumption and the assumptions on $H$ are different from the known results.
We assume \\
($V$) $V\in C^1(\mathbb{R}^3,\mathbb{R})$, there exists $b>0$ such that
\[
 V^b:=\{x\in\mathbb{R}^3:V(x)\leq b\}
\]
has finite measure and denote by $b_{max}:=\sup\{b,|V^b|<\infty\}$.\\
Without confusion, in this part, denote $ \mathbf H:=L^2(\mathbb{R}^3,\mathbb{C}^4)$ and
\begin{equation}\label{eq-the definition of operator A in DE}
A:=-i\displaystyle\sum^3_{k=1}\alpha_k\partial_k+V(x)\beta,
\end{equation}
 which is a unbounded self-adjoint operator on $\mathbf H$,  then we have the following result.
\begin{lem}\label{Bartsch-Ding-2006}\cite[Lemma 3.1]{Bartsch-Ding-2006}
 If $V$ satisfies condition ($V$) then
\[
 \sigma_e(A)\in \mathbb{R}\setminus (-b_{max}, b_{max}).
\]
\end{lem}
Define $\mathcal{L}_s(\mathbf H,-b_{max},b_{max})$ as above, denote $L_s(\mathbb{C}^4)$ the symmetry linear map from $\mathbb{C}^4$ to $\mathbb{C}^4$, here we regard $\mathbb{C}^4$ as $\mathbb{R}^8$.
Redefine $\mathcal{B}\subset C(\mathbb{R}^3,L_s(\mathbb{C}^4))$ here the set of all bounded matrix functions.
Similarly, it is easy to see for any $B\in\mathcal{B}$, it determines a bounded self-adjoint operator on $\mathbf H$, by
 \begin{equation}\label{eq-multiple operator g}
  z(x)\mapsto B(x)z(x), \;\forall z\in\mathbf H,
 \end{equation}
 without confusion, we still denote this operator by $B$, that is to say
 we have the continuous embedding $\mathcal{B}\hookrightarrow \mathcal{L}_s(\mathbf H)$.
Besides, for any  $B_1,B_2\in \mathcal{B}$, $B_1\leq B_2$(or $B_1<B_2$) means that $B_2(x)-B_1(x)$ is semi-positive( or positive) define in $L_s(\mathbb{C}^4)$ for all $x\in \mathbb{R}^3$.
Assume $H\in C^2(\mathbb{R}^3\times\mathbb{C}^4,\mathbb{R})$ satisfying the following conditions
\begin{list}{}
{\setlength{\topsep}{1ex} \setlength{\itemsep}{0ex}
\setlength{\leftmargin}{2.6em}
 }
\item[$(H_0)$] $\nabla_zH(x,0)\equiv 0$, and  $B_0{ :=}\nabla^2_zH(x,0)\in \mathcal{B}\cap \mathcal{L}_s(\mathbf H,-b_{max},b_{max})$.

 \item[$(H_1)$] There exists a constant $ \delta>0$ such that
\[
 -b_{max}+\delta<\nabla^2_zH(x,z)<b_{max}-\delta,\;\forall  (x,z)\in \mathbb{R}^3\times\mathbb{C}^4.
\]

\item[$(H^\pm_2)$] There exist $B_1,B_2\in \mathcal{B}\cap \mathcal{L}_s(\mathbf H,-b_{max},b_{max})$ with $\pm B_1< \pm B_2$,
$i_A(B_1)=i_A(B_2)$ and $\nu_A(B_1)=\nu_A(B_2)=0$,
such that
\[
\pm \nabla^2_zH(x,z)\geq  \pm B_1(x),\;\forall (x,z)\in \mathbb{R}^3\times\mathbb{C}^4
\]
and
\[
\pm \nabla^2_zH(x,z)\leq \pm B_2(x),\;\forall x\in \mathbb{R}^3\;{\rm and}\; |(x,z)|>K
\]
 for some constant $K>0$.
\end{list}
Then we have the following result.
\begin{thm}\label{thm-asymptotically linear result 2}
 Assume (V) is satisfied,  $H\in C^2(\mathbb{R}^3\times\mathbb{C}^4,\mathbb{R})$ satisfies conditions $(H_0)$, $(H_1)$ and $(H^\pm_2)$, if
\[
 i_A(B_0)>i_A(B_2)(or\; i_A(B_0)+\nu(B_0)<i_A(B_1)),
\]
then ($DE$) has a nontrivial solution.  Further more, if $H$ satisfies more conditions, we will get more results.\\
\noindent  A. If $\nabla^2_zH$ is globally Lipschitz continuous on $z$
 and $\nu_A(B_0)=0$, then ($DE$) has another nontrivial solution different from the above one. \\
\noindent B. If $H$ is even in $z$, then ($DE$) has $i_A(B_0)-i_A(B_2)${\rm (}or $ i_A(B_1)-i_A(B_0)-\nu_A(B_0)${\rm )} pairs of nontrivial solutions.
\end{thm}
Compared to the known results, in \cite{Bartsch-Ding-2006,Ding-Ruf-2008} they assume the nonlinear term to be periodic in variable $x$ or satisfying some positive condition.
In Theorem \ref{thm-asymptotically linear result 2}, briefly speaking, we require  that $\nabla^2_z H$ lies in the gape of $\sigma_{ess}(A)$ such that the index pair is well defined. It is the first attempt to study the existence and multiplicity of solutions for Dirac equations via index theory.
\section{The definition of index pair}\label{section-the definition of index pair}

In this section, we will define the index pair $(i_A(B),\nu_A(B))$ by the method of dual variational, then we will  give the relationship between different definitions and the concept of spectral flow.

Let $\mathbf H$ be an infinite dimensional separable Hilbert space with inner product $(\cdot,\cdot)_\mathbf H$ and norm $\|\cdot\|_\mathbf H$.

 Denote by $\mathcal O(\mathbf H)$ the set of all linear self-adjoint operators on $\mathbf H$. For $A\in \mathcal O(\mathbf H)$, we denote by
$\sigma(A)$ the spectrum of $A$ and $\sigma_e(A)$ the essential spectrum of $A$. We define three subsets of $\mathcal O( \mathbf H)$ as follows
$$\begin{aligned} &\mathcal O^-_e(\mu)=\{A\in \mathcal O(\mathbf H)|\;\sigma_e(A)\cap (-\infty,\mu)=\emptyset\;{\rm and}\;\sigma(A)\cap (-\infty, \mu)\ne \emptyset\},\\
&\mathcal O^+_e(\mu)=\{A\in \mathcal O(\mathbf H)|\;\sigma_e(A)\cap (\mu,+\infty)=\emptyset \;{\rm and}\;\sigma(A)\cap (\mu,+\infty)\ne \emptyset\},\\
&\mathcal O^0_e(a,b)=\{A\in \mathcal O(\mathbf H)|\;\sigma_e(A)\cap(a,b)=\emptyset \;{\rm and}\;\sigma(A)\cap (a, b)\ne \emptyset\}.
\end{aligned}$$
We note that if $\mu=+\infty$ and $A\in \mathcal O^-_e(\mu)$, then $\sigma_e(A)=\emptyset$. If $\sigma_e(A)\ne\emptyset$ and $A\in \mathcal O^-_e(\mu)$ for some $\mu$, then $-\infty<\mu<+\infty$ is a real number. Setting $\lambda^-=\inf(\sigma_e(A))$, we have $-\infty<\lambda^-<+\infty$ is real number and $A\in \mathcal O^-_e(\lambda^-)$. Similarly, if $\mu=-\infty$ and $A\in \mathcal O^+_e(\mu)$, then $\sigma_e(A)=\emptyset$. If $\sigma_e(A)\ne\emptyset$ and $A\in \mathcal O^+_e(\mu)$ for some $\mu$, then $-\infty<\mu<+\infty$ is a real number. Setting $\lambda^+=\sup(\sigma_e(A))$, we have $-\infty<\lambda^+<+\infty$ is real number and $A\in \mathcal O^+_e(\lambda^+)$. If the operator $A$ is fixed and $A\in \mathcal O^-_e(\mu)$ or $A\in \mathcal O^+_e(\mu)$, we always write it in $A\in \mathcal O^-_e(\lambda^-)$ or $A\in \mathcal O^+_e(\lambda^+)$ with $\lambda^{\mp}$ in the above sense. We remind that $\inf\emptyset=+\infty$ and $\sup\emptyset=-\infty$.

Let $A\in \mathcal O(\mathbf H)$ satisfying $\sigma(A)\setminus\sigma_{e}(A)\neq\emptyset$. Now, we consider the following cases:

\noindent{\bf Case 1.} $A\in \mathcal O^-_e(\lambda^-)$, { $\lambda^-=\inf(\sigma_e(A))$.}

\noindent{\bf Case 2.}  $A\in \mathcal O^+_e(\lambda^+)$, { $\lambda^+=\sup(\sigma_e(A))$.}

\noindent{\bf Case 3.}  $A\in \mathcal O^0_e(\lambda_a,\lambda_b)$, $-\infty<\lambda_a<\lambda_b<+\infty$.

Denote $\mathcal{L}_s(\mathbf H)$ the set of all  linear bounded self-adjoint operators on $\mathbf H$. Corresponding to  {\it case 1}, {\it case 2} and {\it case 3},
define $\mathcal{L}^-_s(\mathbf H,\lambda^-)$,  $\mathcal{L}^+_s(\mathbf H,\lambda^+)$ and $\mathcal{L}^0_s(\mathbf H,\lambda_a,\lambda_b)$ three subsets of $\mathcal{L}_s(\mathbf H)$ respectively by
\begin{equation}\label{eq-L-}
 \mathcal{L}^-_s(\mathbf H,\lambda^-)=\{B\in \mathcal{L}_s(\mathbf H), \;B< \lambda^-\cdot I\},
\end{equation}
\begin{equation}\label{eq-L+}
 \mathcal{L}^+_s(\mathbf H,\lambda^+)=\{B\in \mathcal{L}_s(\mathbf H), \;B> \lambda^+\cdot I\},
\end{equation}
and
\begin{equation}\label{eq-L0}
 \mathcal{L}_s(\mathbf H,\lambda_a,\lambda_b)=\{B\in \mathcal{L}_s(\mathbf H), \;\lambda_a\cdot I<B< \lambda_b\cdot I\},
\end{equation}
where $I$ is the identity map on $\mathbf H$, $B< \lambda^-\cdot I$ means that there exists $\delta>0$ such that $(\lambda^--\delta)\cdot I-B$ is positive define,
 $B> \lambda^+\cdot I$ and $\lambda_a\cdot I<B<\lambda_b\cdot I$ have  similar meanings.
It is easy to see $\mathcal{L}^-_s(\mathbf H,\lambda^-)$,  $\mathcal{L}^+_s(\mathbf H,\lambda^+)$ and $\mathcal{L}_s(\mathbf H,\lambda_a,\lambda_b)$ are open and convex
subsets of $\mathcal{L}_s(\mathbf H)$.  In \cite{Wang-Liu-2016}, we have defined the index pairs { ($i^{\mp}_A(B),\nu^{\mp}_A(B)$)  and ($i^0_A(B),\nu^0_A(B)$) } in three cases.
In this part, we give a new definition of index pair ($i_A(B),\nu_A(B)$) in the third case.

Now we assume $A\in \mathcal O^0_e(\lambda_a,\lambda_b)$.
Denote by $\mathcal{L}_s(\mathbf H)$ the set of all  linear bounded self-adjoint operators on $\mathbf H$. For any $a,b\in\mathbb{R}$ with $a<b$,
we recall that  $\mathcal{L}_s(\mathbf H,a,b)$ the subset of $\mathcal{L}_s(\mathbf H)$ defined by
\[
 \mathcal{L}_s(\mathbf H,a,b)=\{B\in \mathcal{L}_s(\mathbf H), \;a\cdot I<B< b\cdot I\},
\]
where  $a\cdot I<B< b\cdot I$ means that there exists $\delta>0$ such that $B-(a+\delta)\cdot I$ and $(b-\delta)\cdot I-B$ are positive define.
It is easy to see $\mathcal{L}_s(\mathbf H,a,b)$ is nonempty open and convex subset of $\mathcal{L}_s(\mathbf H)$.
For $A\in \mathcal O^0_e(\lambda_a,\lambda_b)$, we will define the index pairs  ($i_A(B),\nu_A(B)$) for $B\in \mathcal{L}_s(\mathbf H,\lambda_a,\lambda_b)$. Firstly, without any difficulty, we have the following  result as \cite[Lemma 2.1]{Wang-Liu-2016} and we will not prove it here.
\begin{lem}\label{lem-Fredholm property of A-B}  If $A\in \mathcal O^0_e(\lambda_a,\lambda_b)$,
for any  $B\in\mathcal{L}_s(\mathbf H,\lambda_a,\lambda_b)$,
$
 \dim\ker(A-B)<\infty.
$
Further more, if $0\in \sigma(A-B)$, then $0$ is isolated in the point spectrum.
\end{lem}
Now for any $B\in \mathcal{L}_s(\mathbf H,\lambda_a,\lambda_b)$, let $k\in\mathbb{R}\setminus\sigma(A)$ satisfying
\begin{equation}\label{eq-condition of k}
\lambda_a\leq k \;{\rm and}\; k\cdot I<B.
\end{equation}
Consider the bounded self-adjoint operator $ T_{B,k}$ on $\mathbf H$ defined by
\begin{equation}\label{eq-the definotion of operator T}
T_{B,k}:=B^{-1}_k-A^{-1}_k,\;\forall B\in \mathcal{L}_s(\mathbf H,\lambda_a,\lambda_b),
\end{equation}
where \[
B_k:=B-k\cdot I\;{\rm and}\; A_k:=A-k\cdot I.
\]
Firstly, the invertible map $B^{-1}_k$ establishes the one-to-one correspondence between $\ker(T_{B,k})$ and $\ker(A-B)$, so we have $\dim \ker(T_{B,k})=\dim\ker(A-B)$.
Secondly, let $E(z)$ the spectral measure of $A$ and define
\[
P_0:=\displaystyle\int^{\lambda_b-\delta}_{\lambda_a+\delta}1 dE(z),
\]
and
\[
P_1:=I-P_0,
\]
with $\delta$ satisfying $(\lambda_a+\delta)\cdot I<k\cdot I<B<(\lambda_b-\delta)\cdot I$.
Let
\begin{equation}\label{eq-decomposition of H}
\mathbf H=\mathbf H_0\oplus \mathbf H_1,
\end{equation}
with $\mathbf H_*=P_* \mathbf H$($*=0,1$).
It is easy to see
\[
A^{-1}_k|_{\mathbf H_1}<(\lambda_b-k)^{-1}\cdot I,
\]
and
\[
B^{-1}_k>(\lambda_b-k)^{-1}\cdot I.
\]
 we have
\[
(T_{B,k}y, y)_\mathbf H>c(y, y)_\mathbf H,\;\forall y\in \mathbf H_1,
\]
for some fixed $c>0$. Since  $\dim H_0<\infty$,
$T_{B,k}$ has only finite dimensional negative definite subspace, that is to say $(-\infty, 0)\cap \sigma(T_{B,k})$ has only finite points with finite dimensional eigenvalue space.
Summed up, we have the following lemma.
\begin{lem}\label{lem-Decomposition of H}  Suppose  $A\in \mathcal O_e^0(\lambda_a,\lambda_b)$.
For any $B\in \mathcal{L}_s(\mathbf H,\lambda_a,\lambda_b)$ and $k\in\mathbb{R}$ satisfying \eqref{eq-condition of k}, there is an orthogonal decomposition of $\mathbf H$ by
\[
\mathbf H=\mathbf H^-_{T_{B,k}}\oplus \mathbf H^0_{T_{B,k}}\oplus \mathbf H^+_{T_{B,k}},
\]
such that $T_{B,k}$ is negative  definite, zero and positive definite on $\mathbf H^-_{T_{B,k}}$, $\mathbf H^0_{T_{B,k}}$ and $\mathbf H^+_{T_{B,k}}$ respectively. Further more
\[
\dim \mathbf H^-_{T_{B,k}}<\infty,\; \dim \mathbf H^0_{T_{B,k}}=\dim\ker(A-B).
\]
\end{lem}

Thus if $A\in \mathcal O_e^0(\lambda_a,\lambda_b)$, for any $B\in \mathcal{L}_s(\mathbf H,\lambda_a,\lambda_b)$ and $k\in\mathbb{R}$ satisfying \eqref{eq-condition of k}
we denote  the Morse index pair of $T_{B,k}$ by  ($m_A(B),\nu_A(B)$), that is
 \begin{equation}\label{eq-relative Morse index pair}
m_A(B):=\dim \mathbf H^-_{T_{B,k}},\; \nu_A(B):=\dim \mathbf H^0_{T_{B,k}}.
\end{equation}

Of course the index $m_A(B)$ depends on the choose of $k$. But we will show that  $m_A(B_1)-m_A(B_2)$
 will not depend on $k$ for any fixed $B_1,B_2\in  \mathcal{L}_s(\mathbf H,\lambda_a,\lambda_b)$, it only depends on $B_1$, $B_2$ and $A$.
For this purpose, we need the following lemma.
\begin{lem}\label{lem-key lem}
Suppose  $A\in \mathcal O_e^0(\lambda_a,\lambda_b)$. For any $B_1,B_2\in  \mathcal{L}_s(\mathbf H,\lambda_a,\lambda_b)$ satisfying $B_1<B_2$, we have
\[
  m_A(B_2)-m_A(B_1)=\sum_{s\in[0,1)}\nu_A((1-s)B_1+sB_2)
\]
for any $k\in\mathbb{R}\setminus\sigma(A)$ satisfying \eqref{eq-condition of k}.
\end{lem}
\noindent{\bf Proof.}
Denote $i(s):=m_T(B(s))$ and $\nu(s)=\nu(A-B(s))$, where $B(s):=(1-s)B_1+sB_2$. Since $B_1<B_2$, we have $B(s_1)<B(s_2)$, for any $0\leq s_1<s_2\leq1$,
so we have
\[
B^{-1}_k(s_1)>B^{-1}_k(s_2)>0,\; 0\leq s_1<s_2\leq1,
\]
and
\[
 T(s_1)>T(s_2), \; 0\leq s_1<s_2\leq 1,
\]
where $T(s):=B^{-1}_k(s)-A^{-1}_k$ and the map $T(s):[0,1]\to \mathcal{L}_s(\mathbf H)$ is continuous.
Firstly, from the definition of $m_A(\cdot)$, it's easy to see $i(s)$ is left continuous and
\[
 0\leq i(s_1)\leq i(s_2)\leq m_A(B_2),\;\forall\; 0\leq s_1<s_2\leq 1.
\]
Further more, for $s_0\in[0,1]$, if $\nu(s_0)=0$ then $i(s)$ is continuous at $s_0$. If $\nu(s)\neq 0$, we have
$i(s+0)-i(s)=\nu(s)$. In fact, by the continuous  of the eigenvalue of continuous operator function, we have  $i(s+0)-i(s)\le\nu(s)$.
On the other side,  since $T(s_1)>T(s_2)$, for $s_1<s_2$, we see that $i(s+0)-i(s)\ge\nu(s)$.
From the above properties of $i(s)$ and the fact that $i(s)\in[0,m_A(B_2)]\cap\mathbb{Z}$, thus there are only finite number of $s\in[0,1]$ such that $\nu(s)\neq 0$ and
\[
m_A(B_2)-m_A(B_1)=\sum_{s\in[0,1)}\nu(A-(1-s)B_1-sB_2)=\sum_{s\in[0,1)}\nu_A((1-s)B_1+sB_2).
\]
Thus we have proved the lemma.\endproof

Let $\overline{B}:=\displaystyle\frac{\lambda_a+\lambda_b}{2}\cdot I$, then we can define the index pair $( i_A(B),  \nu_A(B))$.
\begin{defi}\label{defi-Definition of index for case 3}
  If $A\in \mathcal O_e^0(\lambda_a,\lambda_b)$, for any $B\in  \mathcal{L}_s(\mathbf H,\lambda_a,\lambda_b)$, define the index pair $( i_A(B),  \nu_A(B))$ by
\begin{align*}
 i_A(B)&:=m_A(B)-m_A(\overline{B}),\\
  \nu_A(B)&:=\dim\ker(A-B).
\end{align*}
\end{defi}
The definition is well defined, we will prove that it only depends on the choice of $\overline B$.
By Lemma \ref{lem-key lem}, for any $\tilde{k}\in\mathbb{R}$ satisfying $ B, \overline{B}<\tilde{k}\cdot I$ and $\tilde{k}<\lambda_b$,
\begin{align*}
m_A(B)-m_A(\overline{B})&=(m_A(\tilde{k}\cdot I)-m_A(\overline{B}))-(m_A(\tilde{k}\cdot I)-m_A(B))\\
               &=\displaystyle\sum_{s\in[0,1)}\nu(A-(1-s)\overline{B}-s\tilde{k}\cdot I)-\sum_{s\in[0,1)}\nu(A-(1-s)B-s\tilde{k}\cdot I),
\end{align*}
where the right hand side does not depend on the choice of $k$ and we have proved that the definition of $i_A(B)$ is well defined.
In this definition, for the fixed operators $\overline{B}$, we have $i_A(\overline{B})=0$. For any other choice of the operators $\overline{B}$, the corresponding index is different up to a constant.

From the definition of the index pair, we can get the following properties.
 \begin{prop}\label{prop-some property of index}
 (1). For any $B\in \mathcal{L}_s(\mathbf H,\lambda_a, \lambda_b)$, $-B\in\mathcal{L}_s(\mathbf H,-\lambda_b, -\lambda_a)$, the index pair {\rm($i_{-A}(-B),\nu_{-A}(-B)$)}  are well defined,
 we have the following equality
\begin{equation}\label{q11} i_A(B)+i_{-A}(-B)+\nu_A(B)\equiv \nu_A(\overline{B}),\;\forall B\in \mathcal{L}_s(\mathbf H,\lambda_a, \lambda_b),\end{equation}
  where $\overline{B}=\frac{\lambda_a+\lambda_b}{2}\cdot I$.

  \noindent(2). For any $B\in \mathcal{L}_s(\mathbf H,\lambda_a, \lambda_b)$, define
  \[
  A_s:=A-(1-s)\overline{B}-sB, \;s\in[0,1].
  \]
  Since $(1-s)\overline{B}+sB\in \mathcal{L}_s(\mathbf H,\lambda_a, \lambda_b)$ for all $s\in[0,1]$, and from Lemma \ref{lem-Fredholm property of A-B}, we have $\{A_s|s\in[0,1]\}$ be a continuous path of self-adjoint Fredholm operators on $\mathbf H$.
The relationship between our index and spectral flows $sf\{A_s,[0,1]\}$ is
\begin{equation}\label{eq12}
sf\{A_s,[0,1]\}=-i_A(B)
\end{equation}
(3). The relationship between the index $i_A(B)$ defined above and  the  index $i^0_A(B)$ defined in \cite{Wang-Liu-2016} is that
\begin{equation}\label{eq13}
i_A(B)+n_0=i^0_A(B),\;\;\forall B\in \mathcal{L}_s(\mathbf H,\lambda_a, \lambda_b),
\end{equation}
where $n_0=\sum_{s\in(0,1)}\nu_A(s\lambda_a+(1-s)\frac{\lambda_a+\lambda_b}{2})$

 \end{prop}
 \noindent{\bf Proof.}(1). Choose $B'\in\mathcal{L}_s(\mathbf H,\lambda_a, \lambda_b)$ satisfying $\lambda\cdot I <B'<B$ and $B'<\overline B$, then we have
 \[
 i_A(B)-i_A(B')=\nu_A(B')+\sum_{\lambda\in(0,1)}\nu_A(\lambda B'+(1-\lambda) B))
 \]
 and
 \[
  i_{-A}(-B')-i_{-A}(-B)=\nu_{-A}(-B)+\sum_{\lambda\in(0,1)}\nu_{-A}(\lambda (-B)+(1-\lambda) (-B))).
 \]
Since $\nu_A(B)=\nu_{-A}(-B)=\dim \ker(A-B)$, we have
\[
i_A(B)+i_{-A}(-B)+\nu_{A}(B)=i_A(B')+ i_{-A}(-B')+\nu_A(B').
\]
Similarly, we have
\[
i_A(\overline{B})+i_{-A}(-\overline{B})+\nu_{A}(\overline{B})=i_A(B')+ i_{-A}(-B')+\nu_A(B').
\]
From the definition, $i_A(\overline{B})=i_{-A}(-\overline{B})=0$, so we have the equality \eqref{q11}.\\
(2). The spectral flow of $A_s$ represents the net change in the number of negative eigenvalue of $A_s$ as $s$ tuns from 0 to 1. If $B>\overline{B}$, by the definition of spectral flow, we have
\begin{equation}\label{eq-prop-1}
sf\{A_s,[0,1]\}=-\sum_{s\in[0,1)}\dim\ker(A_s),
\end{equation}
from Lemma \ref{lem-key lem} and Definition \ref{defi-Definition of index for case 3}, we have $sf(A_s)=-i_A(B)$ for $B>\overline{B}$.
Now, for any $B\in\mathcal{L}_s(\mathbf H,\lambda_a, \lambda_b)$, choose $B'\in\mathcal{L}_s(\mathbf H,\lambda_a, \lambda_b)$ satisfying $B'>\overline{B}$ and $B'>B$. Now, we have three continuous path of self-adjoint Fredholm operators $A_{*,s}$ on $\mathbf H$($*=1,2,3$), with $A_{1,s}=A-(1-s)\overline{B}-sB$, $A_{2,s}=A-(1-s)B-sB'$ and $A_{3,s}=A-(1-s)\overline{B}-sB'$. From the homotopy invariance of spectral flow, we have
\begin{equation}\label{eq-prop-2}
sf\{A_{3,s},[0,1]\}=sf\{A_{1,s},[0,1]\}+sf(\{A_{2,s},[0,1]\}.
\end{equation}
Since $B'>\overline{B}$ and $B'>{B}$ we have
\begin{align}
sf\{A_{2,s},[0,1]\}&=-\sum_{s\in[0,1)}\dim\ker(A_{2,s})\nonumber\\
           &=i_A(B)-i_A(B'),
\end{align}
and $sf\{A_{3,s},[0,1]\}=-i_A(B')$, so we have
\[
sf\{A_{1,s},[0,1]\}=-i_A(B),\;\forall B\in\mathcal{L}_s(\mathbf H,\lambda_a, \lambda_b).
\]
The equality \eqref{eq12} is proved.

(3). The proof of \eqref{eq13} is similar  since Lemma \ref{lem-key lem} is satisfied both for $i_A(B)$ and $i^0_A(B)$. We only need to show that $i_A(\overline{B})=0$, $i^0_A(\overline{B})=n_0$ and these two equalities is from the definitions, we omit the details  here.

\begin{rem}\label{r21}
(1). From the definition, the index pair $(i_A(B),\nu_A(B))$ will also be well defined if $\sigma(A)\cap(\lambda_a,\lambda_b)=\emptyset$. But in this trivial case, for all $B\in\mathcal{L}_s(\mathbf H,\lambda_a, \lambda_b)$, we can prove that $i_A(B)=\nu_A(B)=0$.

\noindent (2).  We can redefine the index $i_A(B)$ by an abstract method via the following three points.
\begin{enumerate}
 \item{}  Define $i_A(\frac{\lambda_a+\lambda_b}{2}\cdot I)=0$.

\item{}    For any $k\in (\lambda_a,\frac{\lambda_a+\lambda_b}{2})$, define
\[
i_A(k\cdot I)=-\dim E[k,\frac{\lambda_a+\lambda_b}{2}),
\]
where $E(z)$ is the spectral measure of $A$, $E[k,\frac{\lambda_a+\lambda_b}{2})$ is the projection map.

\item{}  For any $B\in \mathcal{L}_s(\mathbf H,\lambda_a,\lambda_b)$, define
\[
i_A(B)=\sum_{s\in[0.1)}\nu_A((1-s)k\cdot I+sB)+i_A(k\cdot I),
\]
where $k\in(\lambda_a,\frac{\lambda_a+\lambda_b}{2})$ satisfying $k\cdot I<B$.
\end{enumerate}

It's easy to  prove that this definition happens to coincide with Definition \ref{defi-Definition of index for case 3}.

\noindent (3). Inspired by the idea of relative Morse index (see\cite{Abb-2001,Chang-Liu-Liu-1997,Fei-1995,Long-Zhu-1999}), we can redefined our index by this concept and we will realize it in our follow-up work.
\end{rem}

\section{Proof of Theorem \ref{thm-asymptotically linear result 1} and Theorem \ref{thm-asymptotically linear result 2}}\label{section-Proof of main results}
\subsection{Proof of Theorem \ref{thm-asymptotically linear result 1}}\label{subsection-Homoclinic orbits}
Consider the homoclinic orbit problem of first order nonlinear Hamiltonian system
\[
\left\{\begin{array}{ll}
 \dot{z}(t)=\mathcal{J}\nabla_zH(t,z),\\
\displaystyle\lim_{t\to\infty}z(t)=0,
  \end{array}
\right.\eqno{(HS)}
\]
as mentioned  in Section \ref{section-introduction} and the Hamiltonian function $H$ has the form \eqref{eq-the form of Hamiltonial function H}.
Recall that $A:=-(\mathcal{J}\frac{d}{dt}+L(t))$, then $A$ is self-adjoint on $\mathbf H:=L^2(\mathbb{R},\mathbb{R}^{2N})$
with domain $D(A)=H^1(\mathbb{R},\mathbb{R}^{2N})$ if $L(t)$ is bounded
and  domain $D(A)\subset H^1(\mathbb{R},\mathbb{R}^{2N})$ if $L(t)$ is unbounded.  Before the proof of Theorem \ref{thm-asymptotically linear result 1}, we need the following lemma.

\begin{lem}\label{Lemma 2.3-Ding-Jeanjean-2007}\cite[Lemma 2.3]{Ding-Jeanjean-2007}, $E:=D(|A|^{1/2})$ embeds continuously into $H^{1/2}(\mathbb{R},\mathbb{R}^{2N})$, hence, $E$ embeds continuously into $L^p$ for all $p\geq 2$ and compactly into $L^p_{loc}$ for all $p\geq 1$.
\end{lem}

 Define
\[
 F(z):=\int_{\mathbb{R}}R(t,z(t))dt,\; z\in \mathbf H.
\]
If $R\in C^2(\mathbb{R}\times\mathbb{R}^{2N},\mathbb{R})$ satisfies $(R_1)$, we have $F\in C^1(\mathbf H,\mathbb{R})$. The solutions of  operator equation
\[
 Az=F'(z),\;z\in D(A)\eqno(OE)
\]
are the homoclinic orbits of  (HS). Firstly we consider case ($R^+_2$) in the condition ($R^{\pm}_2$). Choose $\varepsilon>0$ small enough,
such that
\[
B_\varepsilon:=B_1-\varepsilon I\in \mathcal{B}\cap \mathcal{L}_s(L^2,-b_{max},b_{max}),\;i_A(B_\varepsilon)=i_A(B_1),\;\nu(B_\varepsilon)=0\]
 and
\begin{equation}\label{eq-positive condition}
 (B_1-B_\varepsilon)^{-1}-(A-B_\varepsilon)^{-1}=\varepsilon^{-1}\cdot I-(A-B_\varepsilon)^{-1}>0,
\end{equation}
that is to say the operator $ (B_1-B_\varepsilon)^{-1}-(A-B_\varepsilon)^{-1}$ is positive define.
 Define
\[
 A_\varepsilon:=A-B_\varepsilon
\]
and
\[
 F_\varepsilon(z):=F(z)-\frac{1}{2}(B_\varepsilon z,z)_{\mathbf H},\;\forall z\in \mathbf H.
\]
Then we have that the operator equation (OE) is equivalent to
\[
 A_\varepsilon z=F'_\varepsilon(z),\;z\in D(A).
\]
From ($R^+_2$),  $F_\varepsilon\in C^1(\mathbf H,\mathbb{R})$ is convex on $\mathbf H$,
it's Legendre transform $F^*_\varepsilon$ is well defined on $\mathbf H$ and $F^*_\varepsilon\in C^1(\mathbf H,\mathbb{R})$. It is easy to verify that
\[
 F^*_\varepsilon(z)=\int_\mathbb{R}R^*_\varepsilon(t,z(t))dt,\;\forall z\in \mathbf H,
\]
where $R^*_\varepsilon(t,z)$ is the Legendre transform of  $R_\varepsilon(t,z):=R(t,z)-\frac12(B_\varepsilon z,z)$  corresponding to variable $z\in \mathbb{R}^{2N}$.
We have $z\in D(A)$ is solution of (OE) if and only if $u\in \mathbf H$ is a solution of
\[
 A^{-1}_\varepsilon u=F^{*'}_\varepsilon(u),
\]
where $z=A^{-1}_\varepsilon u$. Define the functional $\Psi(u):\mathbf H\to\mathbb{R}$ by
\[
 \Psi(u):=F^*_\varepsilon(u)-\frac{1}{2}(A^{-1}_\varepsilon u,u)_\mathbf H,\;\forall u\in \mathbf H.
\]
So we have the critical point $u$ of $\Psi$ corresponding to the homoclinic orbit $z$ of (HS)
with the relationship given by $z=A^{-1}_ku$.  Let $E:=D(|A|^{1/2})$,
since $\sigma_e(A)\in \mathbb{R}\setminus(-b_{max}+\delta,b_{max}-\delta)$, $0$ is at most an isolate point spectrum of $A$ with finite dimensional eigenspace.
Let $Q:\mathbf H\to \mathbf H$ the projection map on $\ker(A)$, then we can define the norm on $E$ by
\[
 \|u\|^2_E:=(|A|u,u)_\mathbf H+(Qu,u)_\mathbf H,\;\forall u\in E.
\]
\begin{lem}\label{lem-PS condition}
If $R$ satisfies condition ($R_0$), ($R_1$) and ($R^+_2$) then  $\Psi$ satisfies the (PS) condition.
\end{lem}
\noindent{\bf Proof.} Let $\{u_n\}$ be a (PS) sequence of $\Psi$, that is to say $\Psi(u_n)\to c$ and $\Psi'(u_n)\to 0$. We divide the proof into five steps and some steps are from the corresponding part of \cite{Ding-Jeanjean-2007}.

\noindent{\it Step 1.     Show that $\{u_n\}$ is bounded in $\mathbf H$}. Without loss of generality, we can assume $R(t,0)\equiv 0$.  From ($R_1$) and ($R^+_2$),  we have $ F_\varepsilon(z)\leq \frac{1}{2}((B_2-B_\varepsilon)z, z)_\mathbf H+c$ for some constant $c>0$, so
\[
 \Psi(u)\geq \frac{1}{2}((B_2-B_\varepsilon)^{-1}u,u)_\mathbf H-\frac{1}{2}(A^{-1}_\varepsilon u,u)_\mathbf H+c.
\]
Since $(B_1-B_\varepsilon)^{-1}-A^{-1}_\varepsilon>0$, $i_A(B_1)=i_A(B_2)$ and $\nu_A(B_2)=0$,
we have
 \begin{equation}\label{eq-positive condition of operator in PS condition1}
 (B_2-B_\varepsilon)^{-1}-A^{-1}_\varepsilon>0.
 \end{equation}
 If not, there exists $h_0\in \mathbf H\setminus\{0\}$, such that
\[
(((B_2-B_\varepsilon)^{-1}-A^{-1}_\varepsilon)h_0,h_0)_\mathbf H\leq 0.
\]
Let $f(t):[0,1]\to \mathbb{R}$, defined by
\[
f(t)=(((t(B_2-B_1)+\varepsilon)^{-1}-A^{-1}_\varepsilon)h_0,h_0)_\mathbf H,
\]
so we have
\[
f(0)=(((B_1-B_\varepsilon)^{-1}-A^{-1}_\varepsilon)h_0,h_0)_\mathbf H> 0,
\]
and
\[
f(1)=(((B_2-B_\varepsilon)^{-1}-A^{-1}_\varepsilon)h_0,h_0)_\mathbf H\leq 0,
\]
thus, there exists $t_0\in(0,1]$, such that
\[
f(t_0)=(((t_0(B_2-B_1)+\varepsilon)^{-1}-A^{-1}_\varepsilon)h_0,h_0)_\mathbf H=0.
\]
Since $(t_0(B_2-B_1)+\varepsilon)^{-1}-A^{-1}_\varepsilon$ is self-adjoint operator on $\mathbf H$, we have
\[
((t_0(B_2-B_1)+\varepsilon)^{-1}-A^{-1}_\varepsilon)h_0=0
\]
thus we have
\[
Au=(1-t_0)B_1u+t_0B_2u
\]
with $u=A^{-1}_\varepsilon h_0\neq 0$. That is to say $\nu_A((1-t_0)B_1+t_0B_2)\neq 0$ and $i_A(B_2)>i_A(B_1)$ which is contradict to $R^+_2$. Thus we have \eqref{eq-positive condition of operator in PS condition1},
so
\[
 \Psi(u)\to+\infty,\; {\rm as}\;\|u\|_\mathbf H\to\infty.
\]
Since $\Psi(u_n)\to c$, we have $\{u_n\}$ is bounded in $\mathbf H$.
 Denote  $\Phi(z):D(|A|^{1/2})\to \mathbb{R}$ by
\[
\Phi(z):=\frac{1}{2}(Az,z)_\mathbf H-F(z),\;\forall z\in D(|A|^{1/2}),
\]
the variation of (OE) and
\[
z_n:=A^{-1}_\varepsilon u_n,\;\;y_n:=\Psi'(u_n).
\]
\noindent{\it Step 2,  $\{z_n\}$ is bounded in $E(=D(|A|^{1/2}))$ and $\Phi'(z_n)\to 0$ in $E$.}
From the definition of $y_n$, we have
\begin{equation}\label{eq-yn convergent to 0}
y_n\to 0\; {\rm in}\;\mathbf H
\end{equation}
From the definition of $\Psi$ and the property of Legendre transform, we have
\begin{equation}\label{eq-zn yn un}
F^{*'}_\varepsilon(u_n)=z_n+y_n,\;F'_\varepsilon(z_n+y_n)=u_n
\end{equation}
 and
\begin{equation}\label{eq-Ak-Fk1}
 A_\varepsilon z_n=F'_\varepsilon(z_n+y_n),\;{\rm or}\;Az_n=F'(z_n+y_n)-B_\varepsilon y_n.
\end{equation}
From the smoothness of $F'_\varepsilon$ and the relationship between $u_n$, $z_n$ and $y_n$, we have $\{z_n\}$ is bounded in $\mathbf H$.
Since $\nu(B_\varepsilon)=0$, $\mathbf H=\mathbf H^-(A_\varepsilon)\bigoplus \mathbf H^+(A_\varepsilon)$, where $\mathbf H^-(A_\varepsilon)$ and  $\mathbf H^+(A_\varepsilon)$ are the positive define and negative define subspace of $A_\varepsilon$ respectively, for any $z\in \mathbf H$, $z=z^-+z^+$ with $z^*\in \mathbf H^*(A_\varepsilon)$($*=\pm$). From \eqref{eq-Ak-Fk1}, we have
\begin{align*}
\|z_n\|^2_E&=(A_\varepsilon z_n,z_n^+-z_n^-)_\mathbf H\\
           &=(F'_\varepsilon(z_n+y_n),z_n^+-z_n^-)_\mathbf H\\
           &\leq \|u_n\|_\mathbf H\|z_n\|_\mathbf H.
\end{align*}
So $\{z_n\}$ is bounded in $E$. For any $h\in E$, with $\|h\|_E=1$, we have
\begin{align*}
(\Phi'(z_n),h)_E&=(A(z_n)-F'(z_n),h)_\mathbf H\\
              &=(F'(z_n+y_n)-B_\varepsilon y_n-F'(z_n),h)_\mathbf H\\
              &=((\nabla^2_zR(t,z_n+\xi_n)-B_\varepsilon)y_n,h)_\mathbf H\\
              &\leq c\|y_n\|_\mathbf H\|h\|_\mathbf H\\
              &\leq c\|y_n\|_\mathbf H,
\end{align*}
since $y_n\to 0$ in $\mathbf H$ and the arbitrariness of $h$, we have $\Phi'(z_n)\to 0$ in $E$.
The last three steps are the results of \cite[Lemma 4.2-4.4]{Ding-Jeanjean-2007}.

\noindent{\it Step 3. Along a subsequence $\{z_{j_n}\}$, for any $\varepsilon>0$, there exists $r_\varepsilon>0$ such that
\begin{equation}\label{eq-app-1a}
 \displaystyle\limsup_{n\to\infty}\int_{I_n\setminus I_r}|z_{j_n}|^2dt\leq \varepsilon
\end{equation}
for all $r\geq r_{\varepsilon}$.} Note that, for each $n\in\mathbb{N}$, $\int_{I_n}|z_j|^2dt\to\int_{I_n}|z|^2dt $ as $j\to\infty$. There exists $i_n\in\mathbb{N}$ such that
\[
 \displaystyle\int_{I_n}(|z_j|^2-|z|^2)dt<\frac{1}{n},\;\forall  j>i_n.
\]
Without loss of generality, we can assume $i_{n+1}\geq i_n$. In particular, for $j_n=i_n+n$ we have
\[
 \displaystyle\int_{I_n}(|z_{j_n}|^2-|z|^2)dt<\frac{1}{n}.
\]
Observe that there is $r_\varepsilon>0$ satisfying
\begin{equation}\label{eq-app-2b}
 \displaystyle\int_{\mathbb{R}\setminus I_r}|z|^2<\varepsilon
\end{equation}
for all $r\geq r_\varepsilon$. Since
\begin{align*}
 \displaystyle\int_{I_n\setminus I_r}|z_{j_n}|^2&=\int_{I_n}(|z_{j_n}|^2-|z|^2)+\int_{I_n\setminus I_r}|z|^2+\int_{I_r}(|z|^2-|z_{j_n}|^2)\\
                                                                               &\leq \frac{1}{n}+\int_{\mathbb{R}\setminus I_r}|z|^2+\int_{I_r}(|z|^2-|z_{j_n}|^2),
\end{align*}
then we will get \eqref{eq-app-1a}.

Let $\eta:[0,\infty]\to[0,1]$ be a smooth function satisfying
\[
 \eta(s)=\left\{\begin{array}{ll}
                1, &s\leq 1,\\
                0, &s\geq 2.
               \end{array}
\right.
\]
Define $\tilde{z}_n(t)=\eta(2|t|/n)z(t)$ and set $h_n:=z-\tilde{z}_n$.  Since the embedding $E\hookrightarrow L^2_{loc}$ is compact and $z_j(t)\to z(t)$ a.e. in $t$,
we have $z$ is a critical point of $\Phi$. That is to say $z$ is a homoclinic orbit. So we have $h_n\in H^1$ and
\[
 \|h_n\|_E\to 0,\; \|h_n\|_{L^\infty}\to 0,\;{\rm as}\;n\to \infty.
\]
\noindent{\it Step 4.  Show that $\Phi'(z_{j_n}-\tilde{z}_n)\to 0$}. Observe that, $\forall h\in E$,
\[
(\Phi'(z_{j_n}-\tilde{z}_n),h)_E=((\Phi'(z_{j_n})-\Phi'(\tilde{z}_n)),h)_E+\displaystyle\int_\mathbb{R}\nabla_z (R(t,z_{j_n})-R(t,z_{j_n}-\tilde{z}_n)-R(t,\tilde{z}_n))hdt.
\]
Now,  the compactness of Sobolev embeddings implies that, for any $r>0$,
\begin{align*}
 &\displaystyle\lim_{n\to \infty}|\int_{I_r}\nabla_z (R(t,z_{j_n})-R(t,z_{j_n}-\tilde{z}_n)-R(t,\tilde{z}_n))hdt|\\
 \leq& \displaystyle\lim_{n\to \infty}\int_{I_r}|\nabla^2_z(R(t,\xi_n)-R(t,\eta_n))(z_{j_n}-\tilde{z}_n)h|dt\\
 \leq& 2b_{max}\displaystyle\lim_{n\to \infty}\int_{I_r}|(z_{j_n}-\tilde{z}_n)h|dt=0
\end{align*}
uniformly in $\|h\|_E=1$.
 For any $\varepsilon>0$ let $r_\varepsilon>0$ be large enough such that \eqref{eq-app-1a} and \eqref{eq-app-2b} hold.
Then
\[
 \displaystyle\limsup_{n\to \infty}\int_{I_n\setminus I_r}|\tilde{z}_n|^2\leq \int_{\mathbb{R}\setminus I_r}|z|^2\leq \varepsilon
\]
for all $r\geq r_\varepsilon$. From ($R_0$), ($R_1$), \eqref{eq-app-1a} and the fact $\|h\|_{L^2}\leq c\|h\|_E$, we have
\begin{align*}
\displaystyle&\limsup_{n\to\infty}|\int_\mathbb{R}\nabla_z (R(t,z_{j_n})-R(t,z_{j_n}-\tilde{z}_n)-R(t,\tilde{z}_n))hdt|\\
                   =&\limsup_{n\to\infty}|\int_{I_n\setminus I_r}\nabla_z (R(t,z_{j_n})-R(t,z_{j_n}-\tilde{z}_n)-R(t,\tilde{z}_n))hdt|\\
               \leq& c_1\limsup_{n\to\infty}\int_{I_n\setminus I_r}(|z_{j_n}|+|\tilde{z}_n|)|h|dt\\
               \leq& c_1\limsup_{n\to\infty}(\|z_{j_n}\|_{L^2(I_n\setminus I_r)}+\|\tilde{z}_n\|_{L^2(I_n\setminus I_r)})\|h\|_{L^2}\\
               \leq&c_2 \varepsilon^{1/2}.
\end{align*}
Thus we have
\[
 \displaystyle\lim_{n\to\infty}\int_\mathbb{R}\nabla_z (R(t,z_{j_n})-R(t,z_{j_n}-\tilde{z}_n)-R(t,\tilde{z}_n))hdt=0
\]
uniformly for $\|h\|_E=1$ and this proves the result.

\noindent{\it{Step 5. $\{z_j\}$ has a convergent subsequence in $E$}}. Recall the decomposition \eqref{eq-decomposition of H} of $\mathbf H$,  with
\[
 P_0:=\int^{b_{max}-\delta/2}_{-b_{max}+\delta/2}1dE(z),
\]
where $E(z)$ is the spectral measure of $A$, $b_{max}$ and $\delta$ are defined in ($R_1$).
let
\[
 E=E_0\oplus E_1,
\]
with $E_*=E\cap \mathbf H_*$ ($*=0,1$). Let
\[
 x_n:=z_{j_n}-\tilde{z}_n=x_{n,0}+x_{n,1},
\]
with $x_{n,*}\in E_*$ ($*=0,1$).  Thus,  we have
\[
 \|x_{n,1}\|^2_{\mathbf H}\leq \frac{\|x_{n,1}\|^2_E}{b_{max}-\delta/2}.
\]
  Since $x_n\rightharpoonup 0$ in $E$ and $\dim E_0<\infty$, we have $x_{n,0}\to 0$ in $E$, and from the fourth step, $\Phi'(x_n)\to 0$.
Let
\[
 \tilde{x}_{n,1}:=x^+_{n,1}-x^-_{n,1},
\]
where $x^+_{n,1}$ and $x^-_{n,1}$ corresponds to the positive and negative define space of $A$.
From ($R_0$) and ($R_1$), we have
\begin{align*}
 \|\tilde{x}_{n,1}\|^2_E&=\Phi'(x_n)\tilde{x}_{n,1}+\int_{\mathbb{R}}\nabla_z R(t,x_n)\tilde{x}_{n,1}dt\\
                                       &\leq o(1)+\int_{\mathbb{R}}|\nabla^2_z R(t,\xi_n)x_n||\tilde{x}_{n,1}|dt\\
                                       &\leq o(1)+(b_{max}-\delta)\|y_n\|_{\mathbf H}\|\tilde{x}_{n,1}\|_{\mathbf H}\\
                                       &\leq o(1)+\frac{b_{max}-\delta}{b_{max}-\delta/2}\|\tilde{x}_{n,1}\|^2_E.
\end{align*}
Hence we have $\|\tilde{x}_{n,1}\|^2_E\to 0$ and so $\|x_n\|_E\to 0$. Since $z_{j_n}-z=x_n+(\tilde{z}_n-z)$, we have $\|z_{j_n}-z\|_E\to 0$, recall \eqref{eq-yn convergent to 0} and \eqref{eq-zn yn un}, $\{u_n\}$ has a convergent subsequence $\{u_{j_n}\}(u_{j_n}:=F_\varepsilon(z_{j_n}+y_{j_n})$. The proof is complete.\endproof

\noindent{\bf Proof of  Theorem \ref{thm-asymptotically linear result 1}.} \\
\noindent{\it {A. The existence of nontrivial solution.}}\\
 \noindent{\it{Step 1. Consider the case of ($R^+_2$).}} Firstly, without loss of generality, we can assume $R(t,0)\equiv 0$, from the proof of Lemma \ref{lem-PS condition}, we have $\Psi$ is bounded from below, then by Ekeland's variational principle and the (PS) condition, $\Psi$ gets its minimal value at some point $u_0$.
Of  cause, $u_0$ is a critical point of $\Psi$. Secondly, we will prove $u_0 \neq\theta$(the zero element in $\mathbf H$). In fact, from condition $i_A(B_0)>i_A(B_2)$, we have
the operator
\[
(B_0-B_\varepsilon)^{-1}-A^{-1}_\varepsilon
\]
 on $\mathbf H$
has an $i_A(B_0)-i_A(B_2)$- dimensional negative define space, denote it by $Z$.
From condition ($R_0$), we have
\begin{equation}\label{eq-the condition of Psi near 0}
 \Psi(u)\leq \frac{1}{2}((B_0-B_\varepsilon)^{-1}u,u)_\mathbf H-\frac{1}{2}(A^{-1}_\varepsilon u,u)_\mathbf H+o(\|u\|^2_\mathbf H),\; u\in Z\; {\rm and}\; \|u\|_\mathbf H \;{\rm small\; enough}.
\end{equation}
Thus $\theta$ is not a minimal value point of $\Psi$, and $u_0\neq\theta$. \\
\textbf{{\it Step 2. Consider the case of ($R^-_2$).}} In this case we only need to know the following two facts.
One is  that the solutions of operator equation (OE) are also the solutions of the following operator equation
\[
 -Az=-F'(z),\;z\in D(A).
\]
Another is that for $A\in \mathcal{O}^0_e(\lambda_a,\lambda_b)$, if $B\in \mathcal{L}_s(\mathbf H,\lambda_a, \lambda_b)$, we have $-A\in\mathcal{O}^0_e(-\lambda_b,-\lambda_a)$ and
$-B\in\mathcal{L}_s(\mathbf H,-\lambda_b, -\lambda_a)$, from Proposition \ref{prop-some property of index}, we have
\[
 i_A(B)+i_{-A}(-B)+\nu_A(B)=\nu_A(\overline{B}),\;\forall B\in \mathcal{L}_s(\mathbf H,\lambda_a, \lambda_b).
\]
From condition $\nu_A(B_1)=0$ and  $i_A(B_0)+\nu_A(B_0)<i_A(B_1)$, we have  $i_{-A}(-B_0)>i_{-A}(-B_1)$.
So, from the proof of case ($R^+_2$), we can prove this case. \endproof\\
\noindent{{\it B. The existence of another nontrivial solution.}} Since $\nabla^2_zR$ is  globally Lipschitz continuous on $z$, and $F'(z)=\nabla_zR(t,z)$, we have
\begin{align*}
F'(z+h)-F'(z)&=\nabla_zR(t,z+h)-\nabla_zR(t,z)\\
             &=\nabla^2_zR(t,z)h+(\nabla^2_zR(t,z+\xi h)-\nabla^2_zR(t,z))h,\;z,h\in\mathbf H,\xi\in(0,1),
\end{align*}
So we have
\begin{align*}
\|F'(z+h)-F'(z)-\nabla^2_zR(t,z)h\|_\mathbf H&=\|(\nabla^2_zR(t,z+\xi h)-\nabla^2_zR(t,z))h\|_\mathbf H\\
                               &\leq c\|h\|^2_\mathbf H.
\end{align*}
That is to say
\[
F'(z+h)-F'(z)=\nabla^2_zR(t,z(t))h+o(\|h\|_\mathbf H),\;h\to 0.
\] and $F:\mathbf H\to \mathbb{R}$ is $C^2$ continuous, satisfying
\[
(F''(z)x,y)_\mathbf H=\int_\mathbb{R}(\nabla^2_zR(t,z(t))x(t),y(t))dt,\;\forall x,y,z\in\mathbf H.
\]
Similarly, we have $F^*_\varepsilon$ and $\Psi\in C^2(\mathbf H,\mathbb{R})$. Now we only consider the case of ($R^+_2$). From the proof of part A, condition ($R^+_2$) and $\nu_A(B_0)=0$,  we have $\Psi$ is bounded from below, $\theta$ is a non-degenerate critical point with its Morse index $i_A(B_0)-i_A(B_2)>0$,
so by the classical Three-Solution Theorem, we have another critical point which is different from  $\theta$ and $u_0$.\\
\noindent{\it C. The rest part of Theorem \ref{thm-asymptotically linear result 1}}.
In order to prove the rest part of Theorem \ref{thm-asymptotically linear result 1}, we need the following result.
\begin{thm}\label{symmetric mountain lemma 2}\cite{Clark-1972}
Let $\Phi\in C^1(E,\mathbb{R})$ be an even functional on a Banach space $E$. Assume $\Phi(0)=0$ and $\Phi$ satisfies the (PS)-condition, if

\noindent($\Phi_1$) there exists a subspace $E_1\subset E$ with $\dim E_1=j$ and $r>0$ such that
\[
\displaystyle\sup_{u\in E_1\cap S_r}\Phi(u)\leq 0,
\]
($\Phi_2$) there exists a subspace $E_2\subset E$ with ${\rm codim} E_2=k<j$ such that
\[
\displaystyle\inf_{u\in E_2}\Phi(u)>-\infty,
\]
then $\Phi$ has at least $j-k$ pairs of critical points.
\end{thm}

 Let us continuous the proof of Theorem \ref{thm-asymptotically linear result 1}.  If ($R^+_2$) is satisfied, $i_A(B_0)>i_A(B_1)$, $\Psi$ is bounded from below, and satisfies \eqref{eq-the condition of Psi near 0}, with $Z$ defined above an $i_A(B_0)-i_A(B_1)$-dimensional subspace of $\mathbf H$. In Theorem \ref{symmetric mountain lemma 2}, let $E_1=Z$ and $E_2=\mathbf H$, then $\Phi:=\Psi$ satisfies all the conditions, thus we have $\Psi$ has $i_A(B_0)-i_A(B_1)$ pairs of critical points.  If ($R^-_2$) is satisfied, $i_A(B_1)>i_A(B_0)+\nu_A(B_0)$. In this case, let $\Phi:=-\Psi$, with the similar reason, $\Phi$ has $i_{-A}(-B_0)-i_{-A}(-B_1)=i_A(B_1)-i_A(B_0)-\nu_A(B_0)$ pairs of critical points.

\begin{rem}\label{r31}
As mentioned in Remark \ref{r21}, when $\sigma(A)\cap(-b_{max},b_{max})=\emptyset$, we have for any $B_1,B_2\in\mathcal{L}_s(\mathbf H,-b_{max},b_{max})$, $i_A(B_1)=i_A(B_2)=0$. In this case, it says  nothing in Theorem \ref{thm-asymptotically linear result 1}.
Now  we give an example to show
 $\sigma(A)\cap(-b_{max},b_{max})\neq \emptyset$. Let
 \[
 l(t)=
\left\{\begin{array}{ll}
 0,&|t|\leq r_1,\\
 0\leq l(t)\leq b_{max}, &r_1<|t|<r_2,\\
 b_{max},&|t|\geq r_2,
  \end{array}
\right.
\]
 with some $0<r_1<r_2$. Let
 \[
 L(t)=\left(\begin{matrix}
              0 &l(t)\cdot I_N\\
             l(t)\cdot I_N &0
             \end{matrix}\right),\;\;
f(t)=\left\{\begin{array}{ll}
 1,&|t|\leq r,\\
 \frac{r-t}{r_1-r}, &r<|t|<r_1,\\
 0,&|t|\geq r_1,
  \end{array}
\right.
\]
 with some $0<r<r_1$, and $z_0(t)=f(t)\cdot (1,0,0,\cdots,0)\in H^1(\mathbb{R},\mathbb{R}^{2N})$. Then we have
\[
(Az_0,Az_0)_{\mathbf H}=\|z_0'(t)\|^2_{\mathbf H}=2
\]
 and
\[
\|z_0(t)\|^2_{\mathbf H}=2r+\frac{r_1-r}{3}.
\]
So we have $(Az_0,Az_0)_{\mathbf H}<\frac{1}{r}(z_0,z_0)_{\mathbf H}$, that is to say $\sigma(A)\cap [-\frac{1}{\sqrt{r}},\frac{1}{\sqrt{r}}]\neq \emptyset$.
\end{rem}

\subsection{Proof of Theorem \ref{thm-asymptotically linear result 2}}\label{subsection-Dirac equations}
Consider the following Dirac equations
\[
 -i\displaystyle\sum^3_{k=1}\alpha_k\partial_kz+V(x)\beta z=H_z(x,z).\eqno(DE),
\]
introduced in Section \ref{section-introduction}.
Now, in this part let $A:=-i\displaystyle\sum^3_{k=1}\alpha_k\partial_k+V(x)\beta$ which is an unbounded self-adjoint operator in $\mathbf H:=L^2(\mathbb{R}^3,\mathbb{C}^4)$.
Let $E=D(|A|^{\frac{1}{2}})$ be the Hilbert space equipped with th inner product
\[
(u,v)=(|A|^{\frac{1}{2}}u,|A|^{\frac{1}{2}}v)_{L^2}+(u,v)_{L^2}
\]
and norm$\|u\|=(u,u)^{\frac{1}{2}}$. From \cite[Lemma 7.4]{Ding-2007}, we have $E\hookrightarrow H^{\frac{1}{2}}(\mathbb{R}^3,\mathbb{C}^4)$. Let
\[
F(z)=\int_{\mathbb{R}^3}H(x,z(z))dx,
\]
and
\[
F_\varepsilon(z)=\int_{\mathbb{R}^3}H(x,z(z))dx-\frac{1}{2}(B_\varepsilon z,z)_\mathbf H,
\]
where $B_\varepsilon:=B_1-\varepsilon\cdot I$ and $I$ the identity map on $\mathbb{C}^4$. Similar to the above subsection, we can choose $\varepsilon>0$ small enough such that $i_A(B_\varepsilon)=i_A(B_1)$ and $\nu_A(B_\varepsilon)=0$. Let $F^*_\varepsilon$ be  the  Legendre transform of $F$, and $A_\varepsilon:=A-\varepsilon I$, we have $z\in D(A)$ is a solution of ($DE$) if and only if $z$ is a critical point of
\[
\Phi(z):=\frac{1}{2}(A_\varepsilon z,z)_H-F_\varepsilon(z),\;z\in E,
\]
if and only if $u$ is a critical point of
\[
\Psi(z):=\frac{1}{2}(A^{-1}_\varepsilon u,u)_H-F^*_\varepsilon(u),\;u\in \mathbf H,
\]
where $u=A_\varepsilon z$. Similar to Lemma \ref{lem-PS condition}, we have the following result.
\begin{lem}\label{lem-PS condition in DE}
If $F$ satisfies condition $F_0$, $F_1$ and $F^+_2$, then  $\Psi$ satisfies the (PS) condition.
\end{lem}
\noindent{\bf Proof.} Let $\{u_n\}$ be a (PS) sequence of $\Psi$, that is to say $\Psi(u_n)\to c$ and $\Psi'(u_n)\to 0$. Similar to Lemma \ref{lem-PS condition}, we also divide the proof into five steps.
Without any difficult, we can prove the following two steps.\\
\noindent{\it Step 1.} $\{u_n\}$ is bounded in $\mathbf H$.\\
\noindent{\it Step 2.}  $\{z_n\}$ is bounded in $E(=D(|A|^{1/2}))$ and $\Phi'(z_n)\to 0$ in $E$, where
\[
z_n:=A^{-1}_\varepsilon(u_n),
\]
and $\Phi(z):D(|A|^{1/2})\to \mathbb{R}$  is defined by
\[
\Phi(z):=\frac{1}{2}(Az,z)_\mathbf H-F(z),\;\forall z\in D(|A|^{1/2}).
\]
We omit the proof of these two steps here.\\
\noindent{\it Step 3. Along a subsequence $\{z_{j_n}\}$, for any $\varepsilon>0$, there exists $B_\varepsilon>0$ such that
\begin{equation}\label{eq-app-1(DS)}
 \displaystyle\limsup_{n\to\infty}\int_{I_n\setminus I_r}|z_{j_n}|^2dx\leq \varepsilon
\end{equation}
for all $r\geq r_{\varepsilon}$, where $I_r=\{x\in\mathbb{R}^3|\|x\|\leq r\}$.} With the above discussion, we may assume without loss of generality that $z_n\rightharpoonup z$ in $E$. Since  $E\hookrightarrow H^{\frac{1}{2}}(\mathbb{R}^3,\mathbb{C}^4)$, $E\hookrightarrow L^2_{loc}(\mathbb{R}^3,\mathbb{C}^4)$ compactly.  Note that, for each $n\in\mathbb{N}$, $\int_{I_n}|z_j|^2dx\to\int_{I_n}|z|^2dx $ as $j\to\infty$. There exists $i_n\in\mathbb{N}$ such that
\[
 \displaystyle\int_{I_n}(|z_j|^2-|z|^2)dx<\frac{1}{n},\;\forall  j>i_n.
\]
Without loss of generality, we can assume $i_{n+1}\geq i_n$. In particular, for $j_n=i_n+n$ we have
\[
 \displaystyle\int_{I_n}(|z_{j_n}|^2-|z|^2)dx<\frac{1}{n}.
\]
Observe that there is $B_\varepsilon>0$ satisfying
\begin{equation}\label{eq-app-2(DS)}
 \displaystyle\int_{\mathbb{R}^3\setminus I_r}|z|^2<\varepsilon
\end{equation}
for all $r\geq B_\varepsilon$. Since
\begin{align*}
 \displaystyle\int_{I_n\setminus I_r}|z_{j_n}|^2&=\int_{I_n}(|z_{j_n}|^2-|z|^2)+\int_{I_n\setminus I_r}|z|^2+\int_{I_r}(|z|^2-|z_{j_n}|^2)\\
                                                                               &\leq \frac{1}{n}+\int_{\mathbb{R}^3\setminus I_r}|z|^2+\int_{I_r}(|z|^2-|z_{j_n}|^2),
\end{align*}
then we will get \eqref{eq-app-1(DS)}.

Let $\eta:[0,\infty]\to[0,1]$ be a smooth function satisfying
\[
 \eta(s)=\left\{\begin{array}{ll}
                1, &s\leq 1,\\
                0, &s\geq 2.
               \end{array}
\right.
\]
Define $\tilde{z}_n(x)=\eta(2|x|/n)z(x)$ and set $h_n:=z-\tilde{z}_n$.  Since the embedding $E\hookrightarrow L^2_{loc}$ is compact and $z_j(x)\to z(x)$ a.e. in $x$, we have $z$ is a critical point of $\Phi$. That is to say $z$ is a solution of ($DE$) satisfying $z\in W^{1,p}(\mathbb{R}^3,\mathbb{C}^4)(p\geq 2)$. So we have $h_n\in H^1$ and
\[
 \|h_n\|_E\to 0,\; \|h_n\|_{L^\infty}\to 0,\;{\rm as}\;n\to \infty.
\]
\noindent{\it Step 4. We have $\Phi'(z_{j_n}-\tilde{z}_n)\to 0$}. Observe that, $\forall h\in E$,
\[
\Phi'(z_{j_n}-\tilde{z}_n)h=(\Phi'(z_{j_n})-\Phi'(\tilde{z}_n))h+\displaystyle\int_{\mathbb{R}^3}\nabla_z (H(x,z_{j_n})-H(x,z_{j_n}-\tilde{z}_n)-H(x,\tilde{z}_n))hdx.
\]
Now,  the compactness of Sobolev embeddings imply that, for any $r>0$,
\begin{align*}
 &\displaystyle\lim_{n\to \infty}|\int_{I_r}\nabla_z (H(x,z_{j_n})-H(x,z_{j_n}-\tilde{z}_n)-H(x,\tilde{z}_n))hdx|\\
 \leq& \displaystyle\lim_{n\to \infty}\int_{I_r}|\nabla^2_z(H(x,\xi_n)-H(x,\eta_n))(z_{j_n}-\tilde{z}_n)h|dx\\
 \leq& 2b_{max}\displaystyle\lim_{n\to \infty}\int_{I_r}|(z_{j_n}-\tilde{z}_n)h|dx=0
\end{align*}
uniformly in $\|h\|_E=1$.
 For any $\varepsilon>0$ let $B_\varepsilon>0$ be large enough such that \eqref{eq-app-1(DS)} and \eqref{eq-app-2(DS)} hold.
Then
\[
 \displaystyle\limsup_{n\to \infty}\int_{I_n\setminus I_r}|\tilde{z}_n|^2\leq \int_{\mathbb{R}^3\setminus I_r}|z|^2\leq \varepsilon
\]
for all $r\geq B_\varepsilon$. From ($H_0$), ($H_1$), \eqref{eq-app-1(DS)} and the fact $\|h\|_{L^2}\leq c\|h\|_E$, we have
\begin{align*}
\displaystyle&\limsup_{n\to\infty}|\int_{\mathbb{R}^3}\nabla_z (H(x,z_{j_n})-H(x,z_{j_n}-\tilde{z}_n)-H(x,\tilde{z}_n))hdx|\\
                   =&\limsup_{n\to\infty}|\int_{I_n\setminus I_r}\nabla_z (H(x,z_{j_n})-H(x,z_{j_n}-\tilde{z}_n)-H(x,\tilde{z}_n))hdx|\\
               \leq& c_1\limsup_{n\to\infty}\int_{I_n\setminus I_r}(|z_{j_n}|+|\tilde{z}_n|)|h|dx\\
               \leq&c_2 \varepsilon^{1/2}.
\end{align*}
Thus we have
\[
 \displaystyle\lim_{n\to\infty}\int_{\mathbb{R}^3}\nabla_z (H(x,z_{j_n})-H(x,z_{j_n}-\tilde{z}_n)-H(x,\tilde{z}_n))hdx=0
\]
uniformly for $\|h\|_E=1$ and this proves that  $\Phi'(z_{j_n}-\tilde{z}_n)\to 0$.

\noindent{\it Step 5, $\{z_j\}$ has a convergent subsequence in $E$}. Recall the decomposition \eqref{eq-decomposition of H} of $\mathbf H$,  with
\[
 P_0:=\int^{b_{max}-\delta/2}_{-b_{max}+\delta/2}1dE(z),
\]
where $E(z)$ is the spectral measure of $A$, $b_{max}$ and $\delta$ are defined in ($B_1$).
let
\[
 E=E_0\bigoplus E_1,
\]
with $E_*=E\cap \mathbf H_*$ ($*=0,1$). Let
\[
 w_n:=z_{j_n}-\tilde{z}_n=w_{n,0}+w_{n,1},
\]
with $w_{n,*}\in E_*$ ($*=0,1$).  Thus,  we have
\[
 \|w_{n,1}\|^2_{\mathbf  H}\leq \frac{\|w_{n,1}\|^2_E}{b_{max}-\delta/2}.
\]
  Since $w_n\rightharpoonup 0$ in $E$ and $\dim E_0<\infty$, we have $w_{n,0}\to 0$ in $E$, and from above discussion, $\Phi'(w_n)\to 0$.
Let
\[
 \tilde{w}_{n,1}:=w^+_{n,1}-w^-_{n,1},
\]
where $w^+_{n,1}$ and $w^-_{n,1}$ corresponds to the positive and negative define space of $A$.
From ($H_0$) and ($H_1$), we have
\begin{align*}
 \|\tilde{w}_{n,1}\|^2_E&=\Phi'(w_n)\tilde{w}_{n,1}+\int_{\mathbb{R}^3}\nabla_z H(x,w_n)\tilde{w}_{n,1}dx\\
                                       &\leq o(1)+\int_{\mathbb{R}^3}|\nabla^2_z H(x,\xi_n)w_n||\tilde{w}_{n,1}|dx\\
                                       &\leq o(1)+(b_{max}-\delta)\|w_n\|_{\mathbf H}\|\tilde{w}_{n,1}\|_{\mathbf H}\\
                                       &\leq o(1)+\frac{b_{max}-\delta}{b_{max}-\delta/2}\|\tilde{w}_{n,1}\|^2_E.
\end{align*}
Hence we have $\|\tilde{w}_{n,1}\|^2_E\to 0$ and so $\|w_n\|_E\to 0$. Since $z_{j_n}-z=w_n+(\tilde{z}_n-z)$, we have $\|z_{j_n}-z\|_E\to 0$. The proof is complete.\endproof

With the similar argument in the above subsection, we can prove Theorem \ref{thm-asymptotically linear result 2}.

\begin{rem}\label{r32}
Similarly as in Remark \ref{r31}, now
 we give an example to show $\sigma_p(A)\cap (-b_{max},b_{max})\neq\emptyset$. Let $R>0$ large enough, such that $\lambda_1(R)<b^2_{max}$ where $\lambda_1(R)$ is the first eigenvalue of $-\Delta$ on $B_R(0)\subset \mathbb{R}^3$ with Dirichlet boundary condition. Let $z_1$ the corresponding eigenvector and extend it to the whole space by $0$, we still denote it by $z_1$.  Let $V$ satisfying
\[
V(x):=\left\{\begin{array}{ll}
0,&|x|\leq R\\
b_{max},&|x|>2 R
\end{array}
\right.
\]
Then, we have
\begin{align*}
(Az_1,Az_1)_{L^2}&=((-\Delta+V^2+i\sum^3_{k=1}\beta\alpha_k\partial_kV)z_1,z_1)_{L^2}\\
                 &=((-\Delta)z_1,z_1)_{L^2}\\
                 &=\lambda^2_1(R)\|z_1\|^2_{L^2}\\
                 &<b^2_{max}\|z_1\|^2_{L^2}.
\end{align*}
Thus $\sigma(A)\cap (-b_{max},b_{max})\neq \emptyset$.
\end{rem}

 \begin{rem}\label{r33}
 A. As displayed in the proof of our main results, if the functional $\Psi$ is $C^2$ continuous,  we can get more results on the existence and multiplicity of solutions of (HS) and (DE) by Morse theory and critical point theory.\\
 B. The above methods used in subsection \ref{subsection-Homoclinic orbits} and \ref{subsection-Dirac equations} can also be used  to study the existence and multiplicity of solutions of the following diffusion equations
\[
 \left\{
\begin{array}{ll}
\partial_t u-\Delta_x u+V(x)u&=H_v(t,x,u,v),\\
-\partial_t v-\Delta_x v+V(x)v&=H_u(t,x,u,v),
\end{array}
\right.(t,x)\in\mathbb{R}\times \Omega,\eqno(FS)
\]
where $\Omega\subset\mathbb{R}^N$ or $\Omega=\mathbb{R}^N$, $b\in C^1(\mathbb{R}\times\overline{\Omega},\mathbb{R}^N)$, $V\in C(\overline\Omega,\mathbb{R})$ and $H\in C^1(\mathbb{R}\times\overline{\Omega}\times\mathbb{R}^{2N},\mathbb{R})$. (FS) can also be re-written as the form of (OE) and we have the result that if $V$ satisfies some spectral condition then the essential spectrum of $A$ has a gap, thus our methods above can be used here. If $\Omega$ is a bounded domain or $V$ satisfies more stronger condition then  $\sigma_e(A)=\emptyset$ and more theories can be used to study the solution of this situation.$\;$ All of the above will be realized in our subsequence research.\\
C. The condition (L) in subsection \ref{subsection-Homoclinic orbits} is to protect there exists a gap of $\sigma_e(A)$. In an unpublished study, the second author of this paper and Q. Zhang give a weaker condition of $L$ to protect this property, we just give the statement here. Assume \\
($H^\pm_0$)There exist constants $r_0 > 0$ and $b_0 > 0$ such that
\[
\displaystyle\lim_{|s|\to\infty}{\rm meas}\{t\in(s-r_0, s+r_0) | \pm J_0L(t) < b_0\}=0.
\]
Then if $L$ satisfies ($H^+_0$) or ($H^-_0$), we have $\sigma_e(A)\subset\mathbb{R}\setminus(-b_0,b_0)$. Assume ($H^\pm_0$) instead of ($L$), we can also get our results. Similarly, in subsection \ref{subsection-Dirac equations}, we can also give a  weaker condition of $V$ to keep the same result.

All of these will be realized in our following works.
\end{rem}


\begin{thebibliography}{99}

\bibitem{Abb-2001}A. Abbondandolo,
Morse theory for hamiltonian systems,
Chapman \& Hall/CRC, 2001.

\bibitem{Amann-1976} H. Amann,
Fixed point equations and nonlinear eigenvalue problems in ordered Banach spaces,
SIAM Rev. 18 (1976) 620-709.

\bibitem{Amann-Zehnder-1980} H. Amann, E. Zehnder,
Nontrivial solutions for a class of nonresonance problems and applications to nonlinear differential equations,
Annali Scuola Norm. Sup. Pisa  7 (1980) 539-603.

\bibitem{Arioli-Szulkin-1999} G. Arioli, A. Szulkin,
Homoclinic solutions of Hamiltonian systems with symmetry,
J. Differential Equations 158 (1999) 291-313.

\bibitem{Aubin-Ekeland-1984} J.P. Aubin,  I. Ekeland,
 Applied nonlinear analysis,
Wiley, 1984.


\bibitem{Bartsch-Ding-2006}T. Bartsch, Y. Ding,
Solutions of nonlinear Dirac equations,
J. Differential Equations 226 (2006) 210-249.






\bibitem{Chang-1993} K.Q. Chang,
Infinite Dimensional Morse Theory and Multiple Solution Problems,
Birkhauser, Basel, 1993.

\bibitem{Chang-Liu-Liu-1997} K.C. Chang, J.Q. Liu, M.J. Liu,
 Nontrivial periodic solutions for strong resonance Hamiltonian systems,
 Ann. Inst. H. Poincar¨¦ Anal. Non Lin¨¦aire 14 (1997) 103-117.

\bibitem{Chen-Hu-2007} C. Chen, X. Hu,
Maslov index for homoclinic orbits of Hamiltonian systems,
Ann.Inst. H. Poincar\'{e}, Anal. Non lin \'{e}aire  24 (2007) 589-603.

\bibitem{Chen-Ma-2011} G. Chen, S. Ma,
Homoclinic orbits of superlinear Hamiltonian systems,
Proc. Amer.Math.Soc. 139 (2011) 3973-3983.

 \bibitem{Clark-1972}D.C. Clark,
A varant of Ljusternik-Schnirelman theory,
Ind. Univ. Math. J. 22 (1972) 65-74.

\bibitem{Conley-Zehnder-1984}C. Conley,  E. Zehnder,
Morse-type index theory for flows and periodic solutions for Hamiltonian equations,
 Comm. Pure Appl. Math.  37 (1984) 207-253.


 \bibitem{Ding-2006} Y. Ding,
Multiple homoclinics in a Hamiltonian system with asymptotically or super linear terms,
 Commun. Contemp. Math. 8 (2006) 453-480.

\bibitem{Ding-2007} Y. Ding,
Variational Methods for Strongly Indefinite Problems,
World Scientific Publishing, 2007.








 \bibitem{Ding-Jeanjean-2007} Y. Ding, L. Jeanjean,
Homoclinic orbits for a nonperiodic Hamiltonian system,
 J. Differential Equations 237 (2007) 473-490.

\bibitem{Ding-Girardi-1999} Y. Ding, M. Girardi,
Infinitely many homoclinic orbits of a Hamiltonian system with symmetry,
Nonlinear Anal. 38 (1999) 391-415.

\bibitem{Ding-Li-1995} Y. Ding, S. Li,
Homoclinic Orbits for First Order Hamiltonian Systems,
 J. Math. Anal. Appl. 189 (1995) 585-601.


 \bibitem{Ding-Lee-2009} Y. Ding, C. Lee,
Existence and exponential decay of homoclinics in a nonperiodic superquadratic Hamiltonian system,
J. Differential Equations 246 (2009) 2829-2848.

\bibitem{Ding-Liu-2007}Y. Ding, X. Liu,
On Semiclassical Ground States of a Nonlinear Dirac Equation,
Interdiscip. Math. Sci., vol.7, World Scientific, 2007.

\bibitem{Ding-Ruf-2008}Y. Ding, B. Ruf,
Solutions of a nonlinear Dirac equation with external fields,
Arch. Ration. Mech. Anal. 190 (2008) 1007-1032.

\bibitem{Ding-Xu-2015}Y. Ding, T. Xu,
Localized concentration of semi-classical states for nonlinear Dirac equations,
Arch. Ration. Mech. Anal. 216 (2015) 415-447.

\bibitem{Ding-Willem-1999} Y. Ding, M. Willem,
Homoclinic orbits of a Hamiltonian system,
 Z. Angew. Math. Phys. 50 (1999) 759-778,.


\bibitem{Dong-Long-1997}D. Dong, Y. Long,
The iteration formula of Maslov-type index theory with applications to nonlinear Hamiltonian systems,
 Trans. American Math. Soc.  349 (1997) 2619-2661.


\bibitem{Dong-2010} Y. Dong,
Index theory for linear selfadjoint operator equations and nontrivial solutions for asymptotically linear operator equations,
Calc. Var. 38 (2010) 75-109.

\bibitem{Ekeland-1984}I. Ekeland,
Une theorie de Morse pour les systemes hamiltoniens convexes,
Ann IHP Analyse non lineaire 1 (1984) 19-78.


\bibitem{Ekeland-1990} I. Ekeland,
Convexity Methods in Hamiltonian Mechanics,
Springer, 1990.

\bibitem{Ekeland-Hofer-1985}I. Ekeland, H. Hofer,
Periodic solutions with prescribed period for convex autonomous Hamiltonian systems,
 Invent. Math. 81 (1985) 155-188.

\bibitem{Ekeland-Hofer-1987} I. Ekeland, H. Hofer,
Convex Hamiltonian energy surfaces and their closed trajectories,
 Comm. Math. Phys. 113 (1987) 419-467.


\bibitem{Ekeland-Temam-1976} I. Ekeland,  R. Temam,
Convex analysis and variational problems,
 North-Holland-Elsevier, 1976.

\bibitem{Esteban-Sere-1995} M. Esteban, E. S\'{e}r\'{e},
Stationary states of the nonlinear Dirac equation: a variational approach,
Comm. Math. Phys. 171 (1995) 323-350.

\bibitem{Fei-1995} G. Fei,
Relative Morse index and its application to Hamiltonian systems in the Presence of symmetries,
122 (1995) 302-315.

\bibitem{Figueiredo-Pimenta-2017} G.M. Figueiredo, Marcos T.O. Pimenta,
Existence of ground state solutions to Dirac equations with vanishing potentials at infinity,
J. Differential Equations 262 (2017) 486-505.




\bibitem{Hofer-Wysocki-1990} H. Hofer, K. Wysocki,
First order elliptic systems and the existence of homoclinic orbits in Hamiltonian systems,
 Math. Ann. 228 (1990) 483-503.



\bibitem{Liu-2007}C. Liu,
 Maslov-type index theory for symplectic paths with
Lagrangian boundary conditions, Advanced Nonlinear Studies 7 (2007) 131-161.

\bibitem{Liu-2007-2}C. Liu,
  Asymptotically linear Hamiltonian system with Lagrangian
boundary conditions,  Pacific J. Math.  232 (2007) 232-254.

\bibitem{Liu-Long-Zhu-2002} C. Liu, Y. Long, C. Zhu,
Multiplicity of closed characteristics on symmetric convex hypersurfaces in $\mathbb{R}^{2n}$,
 Math. Ann. 323 (2002) 201-215.

\bibitem{Liu-Wang-Lin-2011}  C. Liu, Q. Wang, X. Lin,
 An index theory for symplectic paths associated with two Lagrangian subspaces with applications,
Nonlinearity 24 (2011)  43-70.


\bibitem{Long-1990}Y. Long,
 Maslov-type index, degenerate critical points, and asymptotically linear Hamiltonian systems,
 Sci. China  33 (1990) 1409-1419.

\bibitem{Long-1997}Y. Long,
 A Maslov-type index theory for symplectic paths,
 Topol. Methods Nonlinear Anal.  10 (1997) 47-78.

\bibitem{Long-Zehnder-1990} Y. Long,  E. Zehnder,
Morse theory for forced oscillations of asymptotically linear Hamiltonian systems,
 Stock. Process. Phys. Geom. ed S Alberverio et al(Teaneck, NJ:World Scientific)  (1990) 528-563.

 \bibitem{Long-Zhu-1999} C. Zhu, Y. Long,
Maslov type index theory for symplectiuc paths and spectral flow(I),
 Chinese Ann. of Math. 20B (1999) 413-424.

\bibitem{Long-Zhu-2000-2}Y. Long, C. Zhu,
Maslov type index theory for symplectiuc paths and spectral flow(II),
 Chinese Ann. of Math. 21B (2000) 89-108.

\bibitem{Long-Zhu-2000} Y. Long, C. Zhu,
Closed characteristics on compact convex hypersurfaces in $\mathbb{R}^{2n}$,
 Ann. Math. 155 (2000) 317-368.


\bibitem{Merle-1988} F. Merle,
Existence of stationary states for Dirac equations,
J. Differential Equations 74 (1988) 50-68.


\bibitem{S-1992} E. S\'{e}r\'{e},
Existence of infinitely many homoclinic orbits in Hamiltonian systems,
  Math. Z. 209 (1992) 27-42.

\bibitem{S-1993} E. S\'{e}r\'{e},
Looking for the Bernoulli shift,
 Ann. Inst. H. Poincar\'{e}, Anal. Non lin\'{e}aire 10 (1993) 561-590.

\bibitem{Sun-Chen-Nieto-2011} J. Sun, H. Chen, J. Nieto,
Homoclinic orbits for a class of first-order nonperiodic asymptotically quadratic Hamiltonian systems with spectrum point zero,
 J. Math. Anal. Appl. 378 (2011) 117-127.


\bibitem{Sun-Chu-Feng-2013}J. Sun, J. Chu, Z. Feng,
Homoclinic orbits for first order periodic Hamiltonian systems with spectrum point zero,
Discrete Contin. Dyn. Syst. 33 (2013) 3807-3824.

\bibitem{Szulkin-Zou-2001} A. Szulkin, W. Zou,
Homoclinic orbits for asymptotically linear Hamiltonian systems,
 J. Funct. Anal. 187 (2001) 25-41.

\bibitem{Tanaka-1991} K. Tanaka,
Homoclinic orbits in a first order superquadratic Hamiltonian system: Convergence of subharmonic orbits,
 J. Differential Equations 94 (1991) 315-339.


\bibitem{Wang-Xu-Zhang-2010} J. Wang, J. Xu, F. Zhang,
Homoclinic orbits of superlinear Hamiltonian systems without Ambrosetti-Rabinowitz growth condition,
 Discrete Contin. Dyn. Syst. 27 (2010) 1241-1257.

\bibitem{Wang-Xu-Zhang-2012} J. Wang, J. Xu, F. Zhang,
Infinitely many homoclinic orbits for superlinear Hamiltonian systems,
Topol. Methods Nonlinear Anal. 39 (2012) 1-22.

\bibitem{Wang-Liu-2014}Q. Wang, C. Liu,
Periodic solutions of delay differential systems via Hamiltonian systems,
Nonlinear Ana. TMA 102 (2014) 159-167.

\bibitem{Wang-Liu-2015}Q. Wang, C. Liu,
The relative Morse index theory for infinite dimensional Hamiltonian systems with applications,
J. Math. Anal. Appl.  427 (2015) 17-30.

\bibitem{Wang-Liu-2016} Q. Wang, C.Liu,
A new index theory for linear self-adjoint operator equations and its applications,
J. Differential Equations 260 (2016) 3749-3784


  \bibitem{Zelati-Ekeland-1990} V. Coti Zelati, I. Ekeland, E. S\'{e}r\'{e},
A variational approach to homoclinic orbits in Hamiltonian systems,
 Math. Ann. 228 (1990) 133-160.

\bibitem{Zhang-Tang-Zhang-2013} J. Zhang, X. Tang, W. Zhang,
Homoclinic orbits of nonperiodic superquadratic Hamiltonian system,
Taiwanese J. Math. 17 (2013) 1855-1867.

\bibitem{Zhang-Tang-Zhang-2015} W. Zhang, X. Tang, J. Zhang,
 Homoclinlic solutions for the first-order Hamiltonian system with superquadratic nonlinearity,
 Taiwanese J. Math. 19 (2015) 673-690.






\end{thebibliography}
\end{document}